\documentstyle[twoside,11pt]{article}
\begin{document}
\pagestyle{myheadings}
\markboth{Gantz and  Steer}{On the Riemann-Hilbert-Problem}

\title{Gauge fixing for logarithmic connections over curves
              and the Riemann-Hilbert-Problem\footnote{
   Mathematics Subject Classification (1991) 14H60 (Primary)
      14H30 14F10 14F35 (Secondary)}}

\author{Christian Gantz\thanks{
The first author was fully supported by the
  {\em Alfried Krupp von Bohlen und Halbach -- Stiftung}, Essen}
   \ and Brian Steer \\
   Mathematical Institute, Oxford OX1 3LB, UK \\
   (gantz@maths.ox.ac.uk) }

\maketitle


















\newcommand{\exactss}[5]
{\mbox{$0 \rightarrow #1 \stackrel{#4}{\rightarrow}
#2 \stackrel{#5}{\rightarrow} #3 \rightarrow 0$}}

\newcommand{\fourhorbbb}[4]
{\mbox{$#1 \rightarrow #2 \rightarrow #3 \rightarrow #4$}}

\newcommand{\fourhorbll}[6]
{\mbox{$#1 \rightarrow #2 \stackrel{#5}{\longrightarrow}
#3 \stackrel{#6}{\longrightarrow}
 #4$}}

\newcommand{\fourhorsss}[7]
{\mbox{$#1 \stackrel{#5}{\rightarrow}
#2 \stackrel{#6}{\rightarrow} #3 \stackrel{#7}{\rightarrow}
 #4$}}

\newcommand{\lifta}[6]
{\setlength{\unitlength}{1.0mm}
\begin{center}
\begin{picture}(14,56)(-7,-28)

\putcc{13}{13}{#1}
\putcc{-13}{-13}{#2}
\putcc{13}{-13}{#3}

\multiput(-8,-8)(5,5){3}{\line(1,1){4}}
\put(8,8){\vector(1,1){0}}
\putbr{-1}{-1}{#4}
\put(-6,-13){\vector(1,0){12}}
\puttc{0}{-9}{#6}
\put(13,8){\vector(0,-1){16}}
\putcl{15}{0}{#5}
\end{picture}
\end{center}}

\newcommand{\quadrata}[8]
{\setlength{\unitlength}{1.0mm}
\begin{center}
\begin{picture}(14,46)(-7,-23)

\putcc{-13}{13}{#1}
\putcc{13}{13}{#2}
\putcc{-13}{-13}{#3}
\putcc{13}{-13}{#4}

\put(-6,13){\vector(1,0){12}}
\putbc{0}{15}{#5}
\put(-6,-13){\vector(1,0){12}}
\putbc{0}{-11}{#8}
\put(-13,8){\vector(0,-1){16}}
\putcr{-15}{0}{#6}
\put(13,8){\vector(0,-1){16}}
\putcl{15}{0}{#7}
\end{picture}
\end{center}}

\newcommand{\smallquadrata}[8]
{\setlength{\unitlength}{1.0mm}
\begin{center}
\begin{picture}(14,34)(-7,-17)

\putcc{-13}{8}{#1}
\putcc{13}{8}{#2}
\putcc{-13}{-8}{#3}
\putcc{13}{-8}{#4}

\put(-6,8){\vector(1,0){12}}
\putbc{0}{10}{#5}
\put(-6,-8){\vector(1,0){12}}
\putbc{0}{-6}{#8}
\put(-13,3){\vector(0,-1){6}}
\putcr{-15}{0}{#6}
\put(13,3){\vector(0,-1){6}}
\putcl{15}{0}{#7}
\end{picture}
\end{center}}

\newcommand{\quadratb}[8]
{\setlength{\unitlength}{1.0mm}
\begin{center}
\begin{picture}(14,56)(-7,-28)

\putcc{-18}{13}{#1}
\putcc{18}{13}{#2}
\putcc{-18}{-13}{#3}
\putcc{18}{-13}{#4}

\put(-9,13){\vector(1,0){18}}
\putbc{0}{15}{#5}
\put(-9,-13){\vector(1,0){18}}
\putbc{0}{-11}{#8}
\put(-18,8){\vector(0,-1){16}}
\putcr{-20}{0}{#6}
\put(18,8){\vector(0,-1){16}}
\putcl{20}{0}{#7}
\end{picture}
\end{center}}

\newcommand{\quadratc}[8]
{\setlength{\unitlength}{1.0mm}
\begin{center}
\begin{picture}(14,56)(-7,-28)

\putcc{-23}{13}{#1}
\putcc{23}{13}{#2}
\putcc{-23}{-13}{#3}
\putcc{23}{-13}{#4}

\put(-11,13){\vector(1,0){22}}
\putbc{0}{15}{#5}
\put(-11,-13){\vector(1,0){22}}
\putbc{0}{-11}{#8}
\put(-23,8){\vector(0,-1){16}}
\putcr{-25}{0}{#6}
\put(23,8){\vector(0,-1){16}}
\putcl{25}{0}{#7}
\end{picture}
\end{center}}

\newcommand{\quadratd}[8]
{\setlength{\unitlength}{1.0mm}
\begin{center}
\begin{picture}(14,56)(-7,-28)

\putcc{-29}{13}{#1}
\putcc{29}{13}{#2}
\putcc{-23}{-13}{#3}
\putcc{23}{-13}{#4}

\put(-11,13){\vector(1,0){22}}
\putbc{0}{15}{#5}
\put(-11,-13){\vector(1,0){22}}
\putbc{0}{-11}{#8}
\put(-23,8){\vector(0,-1){16}}
\putcr{-25}{0}{#6}
\put(23,8){\vector(0,-1){16}}
\putcl{25}{0}{#7}
\end{picture}
\end{center}}

\newcommand{\quadrate}[8]
{\setlength{\unitlength}{1.0mm}
\begin{center}
\begin{picture}(14,56)(-7,-28)

\putcc{-25}{13}{#1}
\putcc{29}{13}{#2}
\putcc{-25}{-13}{#3}
\putcc{29}{-13}{#4}

\put(-11,13){\vector(1,0){22}}
\putbc{0}{15}{#5}
\put(-11,-13){\vector(1,0){22}}
\putbc{0}{-11}{#8}
\put(-23,8){\vector(0,-1){16}}
\putcr{-25}{0}{#6}
\put(23,8){\vector(0,-1){16}}
\putcl{25}{0}{#7}
\end{picture}
\end{center}}

\newcommand{\smallquadrate}[8]
{\setlength{\unitlength}{1.0mm}
\begin{center}
\begin{picture}(14,34)(-7,-17)

\putcc{-25}{8}{#1}
\putcc{29}{8}{#2}
\putcc{-25}{-8}{#3}
\putcc{29}{-8}{#4}

\put(-11,8){\vector(1,0){22}}
\putbc{0}{10}{#5}
\put(-11,-8){\vector(1,0){22}}
\putbc{0}{-6}{#8}
\put(-23,3){\vector(0,-1){6}}
\putcr{-25}{0}{#6}
\put(23,3){\vector(0,-1){6}}
\putcl{25}{0}{#7}
\end{picture}
\end{center}}

\newcommand{\quadratf}[8]
{\setlength{\unitlength}{1.0mm}
\begin{center}
\begin{picture}(14,56)(-7,-28)

\putcc{-35}{13}{#1}
\putcc{34}{13}{#2}
\putcc{-35}{-13}{#3}
\putcc{34}{-13}{#4}

\put(-11,13){\vector(1,0){22}}
\putbc{0}{15}{#5}
\put(-11,-13){\vector(1,0){22}}
\putbc{0}{-11}{#8}
\put(-23,8){\vector(0,-1){16}}
\putcr{-25}{0}{#6}
\put(23,8){\vector(0,-1){16}}
\putcl{25}{0}{#7}
\end{picture}
\end{center}}

\newcommand{\threeverta}[5]
{\begin{picture}(140,140)(-175,-70)
\setlength{\unitlength}{0.7pt}
\put(-50,50){$#1$}
\put(-50,-50){$#2$}
\put(-50,-150){$#3$}
\put(-42,37){\vector(0,-1){60}}
\put(-73,0){$#4$}
\put(-42,-63){\vector(0,-1){60}}
\put(65,0){$#5$}
\end{picture}}

\newcommand{\threehorbl}[4]
{\mbox{$#1 \rightarrow #2 \stackrel{#4}{\longrightarrow} #3$}}

\newcommand{\threehorbb}[3]
{\mbox{$#1 \rightarrow #2 \rightarrow #3$}}

\newcommand{\threehorbs}[4]
{\mbox{$#1 \rightarrow #2 \stackrel{#4}{\rightarrow} #3$}}

\newcommand{\threehorsb}[4]
{\mbox{$#1 \stackrel{#4}{\rightarrow} #2 \rightarrow #3$}}

\newcommand{\threehorss}[5]
{\mbox{$#1 \stackrel{#4}{\rightarrow} #2 \stackrel{#5}{\rightarrow} #3$}}

\newcommand{\threehorll}[5]
{\mbox{$#1 \stackrel{#4}{\longrightarrow}
#2 \stackrel{#5}{\longrightarrow} #3$}}

\newcommand{\trianglea}[6]
{\begin{center}
 \begin{picture}(50,48)(-25,-24)
\setlength{\unitlength}{1.0mm}
\put(-18,9){\makebox(0,0){$#1$}}
\put(18,9){\makebox(0,0){$#2$}}
\put(0,-9){\makebox(0,0){$#3$}}
\put(-13,4){\vector(1,-1){8}}
\put(13,4){\vector(-1,-1){8}}
\put(-11,9){\vector(1,0){22}}
\put(0,11){\makebox(0,0)[b]{$#4$}}
\put(-10,-1){\makebox(0,0)[tr]{$#5$}}
\put(10,-1){\makebox(0,0)[tl]{$#6$}}
\end{picture}
\end{center}}

\newcommand{\twohorb}[2]
{\mbox{$#1 \rightarrow #2$}}

\newcommand{\twohors}[3]
{\mbox{$#1 \stackrel{#3}{\rightarrow} #2$}}

\newcommand{\twohorl}[3]
{\mbox{$#1 \stackrel{#3}{\longrightarrow} #2$}}

\newcommand{\N}{\mbox{$I\!\! N$}}	
\newcommand{\Z}{\mbox{$Z\!\!\! Z\!$}}	
\newcommand{\Q}{\mbox{$Q\!\!\!\! I$}}	
\newcommand{\R}{\mbox{$I\!\! R$}}	
\newcommand{\mathR}{R\!\!\!\! I \,\,}

\newcommand{\mathProj}{I\!\! P}
\newcommand{\Proj}{\mbox{$I\!\! P$}}	

\newcommand{\mathC}{C\!\!\!\! I \,\,}
\newcommand{\C}{\mbox{$\, I\!\!\!\! C $}}	

\newcommand{\F}{\mbox{$I\!\! F\!$}}	

\newcommand{\putcc}[3]{\put(#1,#2){\makebox(0,0){$#3$}}}
\newcommand{\putcr}[3]{\put(#1,#2){\makebox(0,0)[r]{$#3$}}}
\newcommand{\putcl}[3]{\put(#1,#2){\makebox(0,0)[l]{$#3$}}}
\newcommand{\puttc}[3]{\put(#1,#2){\makebox(0,0)[t]{$#3$}}}
\newcommand{\puttr}[3]{\put(#1,#2){\makebox(0,0)[tr]{$#3$}}}
\newcommand{\puttl}[3]{\put(#1,#2){\makebox(0,0)[tl]{$#3$}}}
\newcommand{\putbc}[3]{\put(#1,#2){\makebox(0,0)[b]{$#3$}}}
\newcommand{\putbr}[3]{\put(#1,#2){\makebox(0,0)[br]{$#3$}}}
\newcommand{\putbl}[3]{\put(#1,#2){\makebox(0,0)[bl]{$#3$}}}
\newcommand{\la}{\mbox{\bf g}}
\newcommand{\Kahler}{K\"{a}hler}
\newcommand{\calAA}{\mbox{$\cal A$}}
\newcommand{\calE}{\mbox{$\cal E$}}
\newcommand{\calS}{\mbox{$\cal S$}}
\newcommand{\calG}{\mbox{$\cal G$}}
\newcommand{\calK}{\mbox{$\cal K$}}
\newcommand{\calM}{\mbox{$\cal M$}}
\newcommand{\calN}{\mbox{$\cal N$}}
\newcommand{\calL}{\mbox{$\cal L$}}
\newcommand{\calT}{\mbox{$\cal T$}}
\newcommand{\calHH}{\mbox{$\cal H$}}

\newcommand{\calGC}{\mbox{$\cal G$$^{\mathC}$}}
\newcommand{\calNs}{\mbox{$\cal N$$_{S}$}}
\newcommand{\calCs}{\mbox{$\cal C$$_{S}$}}

\newcommand{\calAint}{\mbox{$\cal A$$^{1,1}$}}
\newcommand{\End}{\mbox{\/\/End\/\/}_{0}}
\newcommand{\tcE}{\tilde{\calE}}
\newcommand{\tcF}{\tilde{\calF}}

\newcommand{\tum}{\tilde{U}-V_{U}}
\newcommand{\flag}{\, \mbox{flag} \,}
\newcommand{\delbarF}{\delbar_{\cal{F}}}
\newcommand{\delbarE}{\delbar_{\cal{E}}}


\newcommand{\Aut}{\mbox{Aut}}
\newcommand{\bC}{C\!\!\!\! \underline{I} \,\,}
\newcommand{\blockdiag}{\mbox{block-diag}\,}
\newcommand{\cc}{\mbox{c}}
\newcommand{\calF}{\mbox{$\cal F$}}
\newcommand{\calO}{\mbox{$\cal O$}}
\newcommand{\calW}{\mbox{$\cal W$}}
\newcommand{\diag}{\mbox{diag}\,}
\newcommand{\Endo}{\mbox{End}}
\newcommand{\Gl}{\mbox{Gl}}
\newcommand{\HH}{\mbox{H}}
\newcommand{\Hol}{\mbox{Hol}}
\newcommand{\Image}{\mbox{Im}\,}
\newcommand{\m}[1]{\mbox{#1}}
\newcommand{\mod}{\,\,\, \mbox{mod}\,\,\,}
\newcommand{\normlog}{\mbox{norm log}\, }
\newcommand{\parab}{\mbox{par}\,}
\newcommand{\para}{\mbox{B}\,}
\newcommand{\PD}{\mbox{PD}\,}
\newcommand{\rank}{\mbox{rank}\,}
\newcommand{\Real}{\mbox{Re}\,}
\newcommand{\slope}{\mbox{slope}\,}
\newcommand{\st}{\mbox{ s.t. }\,}
\newcommand{\SU}{\mbox{SU}\,}
\newcommand{\Tr}{\mbox{Tr\,}}
\newcommand{\Tu}{\tilde{U}}
\newcommand{\zP}{z^{\Phi}}
\newcommand{\zmP}{z^{-\Phi}}

\newcommand{\showlabel}[1]{ 
 \label{#1}}
\newcommand{\dimension}{\mbox{dim}}
\newcommand{\stopremark}{ \end{rema} }
\newcommand{\remark}{ \begin{rema} \em  }
\newcommand{\proof}{\noindent{\em Proof:} }
\newcommand{\proofof}[1]{{\noindent \bf Proof} (of #1){\bf :} }
\newcommand{\stopproof}{\\ $\Box$ }
\newcommand{\tu}{\mbox{$\cal{T}(\cal{U})$}}
\newcommand{\ms}{\mbox{$\cal{M}_{S}$}}
\newcommand{\ps}{\mbox{$\cal{P}_{S}$}}
\newcommand{\es}{\mbox{$\cal{E}_{S}$}}
\newcommand{\vv}{\mbox{$\/ \cal V$}}
\newcommand{\uu}{\mbox{$\/ \cal U$}}
\newcommand{\tilu}{\mbox{$\tilde{U}$}}
\newcommand{\M}{\mbox{$\Sigma$}}
\newcommand{\G}{\mbox{$\Gamma$}}
\newcommand{\SS}{\mbox{$\cal S$}}
\newcommand{\D}{\mbox{$\Delta$}}
\newcommand{\k}{\mbox{$\kappa$}}
\newcommand{\h}{\mbox{$\lambda$}}
\newcommand{\g}{\mbox{$\gamma$}}
\newcommand{\Gx}{\mbox{$\Gamma_{x}$}}
\newcommand{\GGx}{\mbox{$G_{x}$}}
\newcommand{\calH}{\mbox{$\cal H$}}

\newcommand{\Gy}{\mbox{$\Gamma_{y}$}}
\newcommand{\Gz}{\mbox{$\Gamma_{z}$}}
\newcommand{\sx}{\mbox{$S_{x}$}}
\newcommand{\sy}{\mbox{$S_{y}$}}
\newcommand{\Sx}{\mbox{$\cal S$$_{x}$}}
\newcommand{\Sy}{\mbox{$\cal S$$_{y}$}}
\newcommand{\Se}{\mbox{$\cal S$$_{0}$}}

\newcommand{\al}{\mbox{$\alpha$}}
\newcommand{\shortfunction}[5]{\mbox{$ #1:#2  \rightarrow #3:%
       #4 \mapsto #5 $} }
\newcommand{\function}[5]{\[ \begin{array}{ccccl}
                              #1 & : & #2 & \rightarrow & #3 \\[1.5mm]
                                 &   & #4 & \mapsto     & #5
                             \end{array}
                          \]}
\newcommand{\action}[4]{\[ \begin{array}{ccl}
                               #1 & \rightarrow & #2 \\[1.5mm]
                               #3 & \mapsto     & #4
                             \end{array}
                          \]}

\newcommand{\del}{\mbox{$\partial$}}
\newcommand{\dd}{\mbox{d}}
\newcommand{\delbar}{\mbox{$\bar{\partial}$}}
\newcommand{\mathdelbar}{\bar{\partial}}
\newcommand{\homo}{\mbox{Hom}}
\newcommand{\olog}{\Omega^{1}( \mbox{log} \,\, 0)}
\newcommand{\logd}{\frac{\dd z}{z}}
\newcommand{\diagmatrix}[2]{\left( \begin{array}{ccccc}
                                    #1 & & & & \\
                                       &\cdot & &0 & \\
                                       & & \cdot & & \\
                                        &0 & & \cdot & \\
                                       & & & & #2
                                   \end{array}
                            \right)
                           }
\newcommand{\uppermatrix}[2]{\left( \begin{array}{cccc}
                                    #1 & & &  \\
                                       &\cdot  & \star  & \\
                                        & 0  & \cdot & \\
                                       & &  & #2
                                   \end{array}
                            \right)
                           }
\newcommand{\twomatrix}[4]{\left( \begin{array}{cc}
                                   #1 & #2 \\
                                    #3 & #4
                                  \end{array}
                            \right)
                            }
\newcommand{\blanctwomatrix}[4]{ \begin{array}{cc}
                                   #1 & #2 \\
                                    #3 & #4
                                  \end{array}
                            }

\newcommand{\bigmatrix}[4]{\left( \begin{array}{cccc}
                                  #1 & \cdots & \cdots & #2 \\
                                  \vdots & \ddots & & \vdots \\
                                  \vdots & & \ddots & \vdots \\
                                   #3 & \cdots & \cdots & #4
                                   \end{array}
                           \right)
                             }
\newcommand{\twovector}[2]{\left( \begin{array}{c}
                                      #1 \\
                                      #2
                                \end{array}
                           \right)
                            }
\newcommand{\threematrix}[9]{\left( \begin{array}{ccc}
                                    #1 & #2 & #3 \\
                                    #4 & #5 & #6 \\
                                    #7 & #8 & #9
                                     \end{array}
                               \right)
                             }

\newcommand{\tz}{\tilde{z}}
\newtheorem{theo}{Theorem}[subsection]
\newtheorem{prop}[theo]{Proposition}
\newtheorem{rema}[theo]{Remark}
\newtheorem{defi}[theo]{Definition}
\newtheorem{coro}[theo]{Corollary}
\newtheorem{lemm}[theo]{Lemma}

\tableofcontents

\newpage


\section{Introduction}

{\sc Said briefly:} We explain in detail the correspondence $\calF$
between algebraic connections over $\Proj^{1}$, logarithmic
at $X =\{ x_{1},...,x_{n} \} \subseteq \Proj^{1}$, and flat bundles
over $\Proj^{1} -X$ with integer weighted filtrations
near each $x_{j}$. Included is a gauge fixing theorem
for logarithmic connections. (Thus far, one could work over any
Riemann surface.)
We prove a bound on the splitting type of a semi-stable
logarithmic connection over $\Proj^{1}$. Using
this we extend some results on the
Riemann-Hilbert-Problem and explain some others.
The work is self contained and elementary, using only basic knowledge
of gauge theory and the Birkhoff-Grothendieck-Theorem.

\hfill

{\sc The concepts:}
 A {\em logarithmic connection} over $(\Proj^{1},X)$ consists of a
holomorphic vector bundle $E \rightarrow \Proj^{1}$ with an algebraic
connection
\[ \nabla: \Omega^{0}(E) \rightarrow \Omega^{0}(E) \otimes
      \Omega^{1}_{\mathProj^{1}}(\log X), \]
satisfying the Leibnitz rule,
where $\Omega^{1}_{\mathProj^{1}}(\log X)$ is the sheaf of
holomorphic 1-forms generated near $x_{j}$ by $\dd z_{j}/z_{j}$
for a coordinate $z_{j}$ centred at $x_{j}$.
$H:=(E,\nabla)|_{\mathProj^{1}-X}$ is a flat bundle.
Isomorphism classes of flat bundles of rank $r$
correspond to conjugacy classes of representations
$\chi:\pi_{1}(\Proj^{1}-X) \rightarrow \Gl \, (r,\C)$,
\cite[p 200]{ati}, \cite[p 51-56]{ano}, \cite[p 4]{kob}.
$\chi$ is called the
monodromy (or holonomy) of $H$.

If $E$ is trivial, one calls $(E,\nabla)$ a
{\em Fuchsian system}.
Choose a global coordinate $z$ on $\Proj^{1}$ such that
$a_{j}:=z(x_{j})\neq \infty$.
For any Fuchsian system $(\Proj^{1} \times \C^{r}, \nabla)$
there exist $B_{j} \in \mbox{ End }(\C^{r})$, \cite[p 4]{ano}, such that
\begin{equation}
 \nabla=\dd + \sum_{j=1}^{n} \frac{B_{j}}{z-a_{j}}.
\label{fz}
\end{equation}

\hfill

{\sc The problem:}
In 1900 Hilbert stated his twenty first problem:
{\em Prove that for any given singularities $X$ and representation
$\chi$ there exists a Fuchsian system realising $(X,\chi)$.}
Literally, \cite{hil},
he said {\em Fuchsian equation}, i.e.\ higher order
differential equations with prescribed singularities. But Anosov \&
Bolibruch argue that he meant
vector-valued linear equations, i.e.\ Fuchsian systems, because the
alternative was already known to be wrong in 1900.
Fuchsian equations induce Fuchsian systems, \cite[Ch. 7]{ano}.
Since Riemann worked on the
problem earlier,
it is called the Riemann-Hilbert-Problem (RHP).
For a comprehensive collection of known results and references
to the RHP
see \cite{ano}, also \cite{bea} and \cite{bol}.
Much of the recent work  is due to Bolibruch.
(An approach different from most is Hain's, \cite{hai}.)
Bolibruch discovered a pair $(X,\chi)$, of rank $r=3$ and
with $n=4$, which cannot be
realised by any Fuchsian system, \cite[p 74-76]{bol}, \cite[p 14]{ano}.
Therefore, he modified the RHP to the question of which
$(X,\chi)$ can occur on Fuchsian systems.
Bolibruch shows that for fixed $\chi$ but varying $X$, the answer
can be different.
We do not consider the dependence on $X$ and
concentrate on positive answers to the RHP.

By the Birkhoff-Grothendieck-Theorem (BGT), \cite{oss},
any vector bundle
$E \rightarrow \Proj^{1}$ is isomorphic to
$\calO(c_{1}) \oplus ... \oplus \calO(c_{r})$ for unique integers
$c_{1} \geq ... \geq c_{r}$, called the {\em splitting type} of $E$.
So, considering the space of all logarithmic connections over
$(\Proj^{1},X)$, the Fuchsian systems (\ref{fz}) constitute
the connected
component of the trivial connection.
Fuchsian systems are clearly
semi-stable.

\hfill

{\sc The approach:}
We follow Deligne, \cite{del}.
To each $x_{j}$, let $U_{j}$ be a small simply-connected
neighbourhood  and $U_{j}^{*}:=
U_{j} - \{ x_{j} \}$.
We show directly that a logarithmic connection $(E,\nabla)$
admits, over a small open neighbourhood of $x_{j}$,
a {\em normal trivialisation}
(Definition \ref{nt}).
This is used to construct on $H=(E,\nabla)|_{\mathProj^{1}-X}$
a filtration
$0 \subset H_{j}^{1} \subset ... \subset H_{j}^{l_{j}}=H|_{U_{j}^{*}}$
by flat subbundles with integer weights
$\Phi_{j}=\mbox{ diag }(\phi_{j}^{i})$,
$(\phi_{j}^{1} \geq ... \geq \phi_{j}^{r}) \in \Z\, ^{r}$, for each
$j=1,...,n$.
Conversely, such data on a flat bundle
$H \rightarrow (\Proj^{1}-X)$ induces
a unique extension of $H$ to a logarithmic connection
$(E,\nabla):= \calF(H,H_{j}^{m},\Phi_{j})$ over $(\Proj^{1},X)$.

Extending and restricting appropriate morphisms,
$\calF$ becomes an equivalence between the categories
of weighted flat bundles $(H,H_{j}^{m},\Phi_{j})$ over $\Proj^{1}-X$
and the category of logarithmic connections
$(E,\nabla) \rightarrow (\Proj^{1},X)$.
The equivalence $\calF$ has been constructed
slightly differently by Manin,
\cite{man};
Deligne, \cite{del}, and Simpson, \cite{sim}, and on
objects partially by Anosov \& Bolibruch, see also \cite{for}.
$\calF$ prerves injections and surjections.
The integer weights are used to define the degree of a weighted flat
bundle.
By Simpson, $\calF$ preserves
degrees and hence (semi-) stability (Definition
\ref{fd}).

If $\gamma_{j}$ is a loop in $U_{j}^{*}$ going once around $x_{j}$, the
parallel transport in $H$ w.r.t. $\gamma_{j}$ is conjugation
equivalent to an upper-triangular matrix. So,
filtrations of $H|_{U_{j}^{*}}$ by flat subbundles exist. There is much
freedom in choosing integer weights.
Hence, any pair $(X,\chi)$ is realized by
several logarithmic connections.
(This even holds over Riemann surfaces, \cite{roe}.)
If one is satisfied with any logarithmic connection
realizing a given pair $(X,\chi)$, the problem is therefore solved;
the difficulty is to decide when the underlying bundle is trivial.

We seek, for given H, filtrations  $H_{j}^{m}$
and integer weights $\Phi_{j}$ such that $\calF(H,H_{j}^{m},\Phi_{j})$
is Fuchsian.
To indicate the relation between our approach and previous ones,
let $(X,\chi)$ be realised by $(E,\nabla)$.
$E$ admits a system $W=(w_{1},...,w_{r})$ of global meromorphic
section, holomorphic away from $x_{1}$,
spanning $E$ off $x_{1}$. $W$ generates a flat bundle over
$\Proj^{1}$. So, every pair $(X,\chi)$ is realized
by a {\em regular system}, i.e.\ a singular algebraic connection
on $\Proj^{1} \times \C^{r}$
such that the flat sections have at most polynomial
growth. This has long been known, \cite{ple}, \cite{del},
and most attempts to find
Fuchsian systems are by ``modifying''
(see \cite[p 77]{ano}) regular ones.
Conversely, regular systems induce logarithmic connections.
To see this,
use the system of sections $V$ as in equation (2.2.21) of \cite{ano}
to generate a free rank $r$ sheaf, i.e.\ vector bundle, and apply
Levelt's result, \cite[p 28]{ano}, \cite[p 379]{lev}.
The modification of regular systems does correspond to changing
filtrations and integer weights on $H$. Bolibruch essentially
introduced the approach, but worked himself mainly via regular systems.

Bolibruch found that any irreducible representation is the monodromy
of a Fuchsian system for any given singularities, \cite[p 83]{ano}.
Having this, one attempts
to apply induction on reducible ones. The difficulty is that the smaller
subspaces in the local filtrations of $H$
have higher integer weights and tend to be
contained in global subspaces. This restricts the choice of filtrations
and weights which make $H$ into a semi-stable weighted flat bundle;
which is neccessary should $H$ be a restriction of a Fuchsian system.
This difficulty comes up in Theorems \ref{re}
and \ref{p} and, in an extreme
form in Proposition \ref{bt}.
The results on reducible representations that we have, Lemma \ref{cc},
follow from
 the preservation of
short exact sequences under $\calF$ and the BGT.

\hfill

{\sc What is new:}
We give a direct proof of a gauge fixing theorem
for logarithmic connections
over curves, Theorem \ref{nf}.
The description of the inverse of  $\calF$ via this gauge fixing theorem
seems new.

We work on the RHP via logarithmic connections, avoiding
regular systems. Instead of Bolibruch's {\em sum of exponents}
of a regular system we use the degree of a bundle
over $\Proj^{1}$ and the concept of semi-stability.
In particular, the preservation of semi-stability under $\calF$ is
usefull because any Fuchsian system is semi-stable.
Bolibruch does not
mention the concept of semi-stability in relation to the RHP.
Applying the properties of $\calF$, explained in the first part of this
article, several of Bolibruch's results on
the RHP follow easily from the Birkhoff-Grothendieck-Theorem and
the fact that $\HH^{0}(\calO(c))$ equals $0$ if
$c<0$ and $\C$ if $c=0$.
We do not reprove this way as many results as possible,
restricting to some signific ones, e.g.
Theorem \ref{tw}, Proposition \ref{bt},  Lemma
\ref{cc} here and \cite[Lem. 5.2.2]{ano} and \cite[Thm. 5.2.2]{ano}.

Besides, perhaps, a more conceptual
proof of known results, we have new ones.
Bolibruch's first counter-example to the RHP
implies that a semi-stable logarithmic connection is not
neccessarily Fuchsian.
However, we prove that any semi-stable logarithmic
connection $(E,\nabla)$ has bounded splitting type, Theorem \ref{ne}.
To be precise,
$c_{i}-c_{i+1} \leq n-2$ for $i=1,...,r-1$ where
$E=\calO(c_{1}) \oplus ... \oplus \calO(c_{r})$,
$c_{1} \geq ... \geq c_{r}$,
$n=\sharp X$.
Bolibruch treats the special case
of logarithmic connections with irreducible
monodromy. His bound, Corollary \ref{bo} here,
 is weaker.
Combined with a technical result of Bolibruch, namely
Proposition \ref{bq} here, i.e.\ \cite[Lem. 4.1.3]{ano}, of which
we provide a direct proof, our bound leads to
the existence of a Fuchsian system with given monodromy
$\chi:\pi_{1}(\Proj^{1}-X) \rightarrow \mbox{Aut}\, (\C^{r})$
if some $\chi(\gamma_{k})$ admits an eigenvector which is a cyclic
vector of the $\pi_{1}(\Proj^{1}-X)$-module $\C^{r}$.
Firstly, this implies Bolibruch's positive solution for irreducible
representations.
Secondly, it gives a shorter proof of his
result that each $\chi$ is a subrepresentation of the monodromy
of a Fuchsian system of double the rank. Thirdly, it leads to an
alternative proof of his complete answer to the RHP in rank three.

We have a new,
sufficient condition for parabolic representations to come from
Fuchsian systems, Theorem \ref{p},
and show that this is always satisfied
in rank four.
A new result for reducible representations
is part (ii) of Theorem \ref{re}.

\hfill

{\sc Acknowledgements:}
We thank Michael Thaddeus for  bringing the RHP to our attention and
showing us Bolibruch's work.

\section{\sloppy
Local logarithmic connections and weighted flat bundles}

\subsection{Logarithmic connections}

Let $U$ be a simply connected, open neighbourhood of $0 \in \C$, $z$
the natural coordinate.
By $\Omega_{U}^{p}$ we denote the sheaves of holomorphic forms on $U$,
Let
\[ \Omega_{U}^{1}(\log 0):=\Omega_{U}^{0} \cdot (\logd) \]
be the free {\em sheaf
of holomorphic 1-forms logarithmic at $0$}, \cite[p 449]{gah}.
Naturally, this sheaf is isomorphic to $\Omega^{0}_{U}(\calK_{U}
\otimes [0])$, where
$\calK_{U}$ is the canonical bundle with section $\dd z$
 and $[0]$ the line bundle on $U$
associated with the divisor $0$.
The vector $(\dd z/z)(0) \in (\calK_{U} \otimes [0])_{0}$
is independent of the choice of coordinate because,
if $u:U \rightarrow \C$ is another coordinate with $u(0)=0$,
\[ \frac{\dd u}{u}=\frac{z \, \dd u}{\dd z} \cdot \frac{\dd z}{z} =(1+o(z))
\frac{\dd z}{z}.\]

\begin{defi}
\showlabel{lc}
A (local) connection logarithmic at 0 is a holomorphic
vector bundle $E \rightarrow U$
and  a $\C$-linear map
\[ \nabla : \Omega^{0} (E) \rightarrow \Omega^{0} (E)
\otimes \Omega_{U}^{1}(\log 0) \]
satisfying the Leibnitz rule $\nabla (f v) =(\dd f ) v + f \nabla (v)$
for all $f \in \Omega^{0}$ and $v \in \Omega^{0} (E)$.
A morphism $\tau: (E',\nabla') \rightarrow (E,\nabla)$
of logarithmic connections is a bundle map
 such that
\[ \tau_{*} \circ \nabla' = \nabla \circ \tau_{*} \]  for
$\tau_{*}:\Omega^{0}(E') \rightarrow \Omega^{0}(E)$.
\end{defi}

A morphism is called {\em injective (surjective)}
if it is as bundle map.
A short sequence
$ \threehorbb{(E',\nabla')}{(E,\nabla)}{(E'',\nabla'')} $
is called {\em exact} if it is as sequence of bundle maps.
Simpson calls $(E,\nabla)$ a regular singular $D_{U}$-module, \cite{sim}.

If $v \in \Omega^{0}(E)$ and $f \in \Omega^{0}$ then
\[ \nabla(fv)=(\frac{\dd f}{\dd z})z v \logd + f\nabla(v)
\mbox{\ \ so \ \ }  (\nabla(fv))(0) = f(0) (\nabla(v))(0) .\]
Hence, $\nabla$ induces a canonical endomorphism
\[  \rho':E_{0} \rightarrow E_{0} \otimes (\calK_{U} \otimes [0])_{0}.  \]

\begin{defi}
The canonical map
$ \rho:E_{0} \rightarrow E_{0}$,
determined by
$ \rho'(w) = \rho(w) \cdot (\dd z/z)(0)$
for all $w \in E_{0}$, is called the residue of $\nabla$ at $0$.
If $(\lambda^{1},...,\lambda^{r})$ denote the eigenvalues of $\rho$
and
$ \phi^{i}:= [- \mbox{Re} (\lambda^{i})] \in \Z$
then
$ \phi^{1} \geq ... \geq \phi^{r} $
 are called the integer weights of $\nabla$ (at 0).
\end{defi}

 We encode the integer weights as
\[ \Phi:= \mbox{diag} (\phi^{i}) = \mbox{block-diag}(\psi^{m} I_{d^{m}}) \]
 with $\psi^{i} > \psi^{i+1}$ for all $i$.
$\Phi$ induces
 a canonical block structure on all matrices which we will use a lot.
If $\nabla_{0}$ is a second logarithmic connection on $E$, the Leibnitz rule
implies that
\[(\nabla_{0} - \nabla ) : E \rightarrow E \otimes \calK_{U} \otimes
     [0] \]
is a
holomorphic bundle map. So, in a
 trivialisation $\theta:E \rightarrow U \times \C^{r}$,
\[ \nabla_{\theta} := \theta \circ \nabla \circ \theta^{-1} =
  \dd + A(z) \frac{\dd z}{z} \]
for holomorphic $A : U \rightarrow \Endo (\C^{r})$.
The converse clearly holds.
Furthermore,
$ \rho=\theta(0)^{-1} \circ A(0) \circ \theta(0)$.

\subsection{Gauge fixing for logarithmic connections}

If all eigenvalues $\mu$ of $K \in \Endo (\C^{r})$ satisfy
$\Real (\mu) \in [0,1)$ then we say
$K$ has {\em normalised eigenvalues}.
If also $G=\exp(2 \pi i K)$ we call $K$ the {\em normalised logarithm}
of $G$: $K =  \mbox{norm log}\,\,  G$.
If $G$ is upper-triangular and has only one eigenvalue
$\rho$ then, for $\mu= \mbox{norm log}\, \rho$, we have, \cite[p 376]{lev},
\[ \mbox{norm log }G=\mu I+\frac{1}{2 \pi i}
  \sum_{j=1}^{\infty} \frac{(-1)^{j}}{j} (\frac{1}{\rho} G -I)^{j}. \]

\begin{defi}
\showlabel{nt}
A trivialisation
$\theta:(E,\nabla) \rightarrow (U \times \C^{r},\nabla_{\theta}
)$ is called normal (w.r.t. $z$) if
\[ \nabla_{\theta} = \dd + z^{\Phi} (-K -\Phi )z^{- \Phi} \logd \]
for some constant, block-upper-triangular $K \in \Endo (\C^{r} )$
with normalised eigenvalues, where $\Phi$ is the integer
weights-matrix of $(E,\nabla)$.
\end{defi}

For a normal trivialisation $\theta$, integer weights
$\Phi=\mbox{ block-diag }(\psi^{m} I_{d^{m}})$ and $(e_{1},...,e_{r})$
the standard frame of $U \times \C^{r}$, let
\[ F^{m}:=\langle e_{d^{1}+...+d^{m-1}+1},...,e_{d^{1}+...+d^{m}}
\rangle \subseteq U \times \C^{r} \]
for $m=1,...,l$ and $F^{0}:= U \times \{ 0 \}$.
Set
\function{\phi}{U\times \C^{r}}{\Z \cup \{ + \infty \}} {v}%
{\left\{ \begin{array}{lcl}
                  \psi^{ m} & \mbox{if} & v\in (\oplus_{0}^{m} F^{k})
                                - (\oplus_{0}^{m-1} F^{k})  \\
                                           + \infty & \mbox{if} & v=0.
                                           \end{array}
                                 \right.}
Clearly, $\phi$ is invariant under parallel
transport away from the singularity.

\remark
\showlabel{ag}
For $v \in U^{*} \times \C^{r}$,
$\phi(v)$ can be described as the integer part of
the {\em asymptotic growth},
\cite[p 17]{ano}, \cite[p 374]{lev}, of the flat extension of $v$ over
$U^{*}$. This follows from Lemma \ref{fs} and equations
(2.2.8), (2.2.11) and (2.2.12) in \cite{ano}. Anosov \&
Bolibruch use Levelt's work on asymptotic growth for
{\em regular systems}, i.e.\ singular connections such that flat sections
have at most polynomial growth. We avoid regular systems.
\stopremark

Under the equivalence of asymptotic growth and $\phi$, (i) and the first part
of (iii) of the following theorem correspond to results
of Levelt, \cite[p 28]{ano}, \cite[p 60]{bol}, \cite[p 379]{lev}.
Gantmacher has
(i), \cite[p 185,191]{gat}.

\begin{theo}
\showlabel{nf}
\begin{description}
\item[(i)]
For each logarithmic connection (and any coordinate),
there exists a normal trivialisation in
some small neighbourhood of the singularity.
\item[(ii)]
Let $\tau:(E',\nabla') \rightarrow (E,\nabla)$ and consider two
normal trivialisations
$\theta:(E,\nabla) \rightarrow (U \times \C^{r},\nabla_{\theta})$ w.r.t.
$z$ and
$\theta':(E',\nabla') \rightarrow (U \times \C^{r'},\nabla'_{\theta'})$
w.r.t. $u$. Then $M:=\theta \circ \tau \circ (\theta')^{-1}$ satisfies
\[ \phi(M(v')) \geq \phi'(v') \] for all $v' \in U \times \C^{r'}$.
\item[(iii)]
If $\tau$ is injective then $\phi(M(v'))=\phi'(v')$ for all $v' \in
U \times \C^{r'}$. If $\tau$ is surjective and $v \in U \times \C^{r}$ then
there exists $v' \in \tau^{-1}(v)$ such that $\phi(v)=\phi'(v')$.
\end{description}
\end{theo}

We give a direct proof.

\proof
{\bf (i):}
Start with any trivialisation and write the connection as
$ \dd + A(z) dz/z$
for $A(z) = \sum_{0}^{\infty} A^{j} z^{j}$.
Applying  a constant gauge transformation,
we can assume
$ A^{0} = \mbox{ block-diag }(A_{m,m}^{0} )$
where each eigenvalue $\lambda$ of $A_{m,m}^{0}$ satisfies
\begin{equation}
 [ - \Real \, \lambda ] = \psi^{m},
\mbox{ \ \ i.e. \ \ }
 - \Real \, \lambda - \psi^{m} \in [0,1).
 \label{three}
\end{equation}
Assume we could find $M(z)=\sum_{0}^{\infty} M^{j} z^{j} :U
\rightarrow \Gl (r,\C)$ and \newline
$B(z)= \sum_{0}^{\infty} B^{j} z^{j} :U
\rightarrow \Endo (\C^{r})$ such that
\begin{equation} \label{one}
 M^{-1} \circ (\dd + A(z) \logd ) \circ M =
\dd + B(z) \logd
\end{equation}
with $M^{0} = I$, hence $B^{0} = A^{0}$, and
$ B(z) = z^{\Phi} (-K -\Phi) z^{-\Phi}$
for $K$ constant, block-upper-triangular.
Then the eigenvalues of $K$
would be those of $-B^{0} - \Phi = -A^{0} - \Phi $
and hence $K$ would have
 normalised eigenvalues by (\ref{three}). We would be done
if the series of $M$ converges in a small neighbourhood of $0$.

Eq. (\ref{one}) is equivalent to
\[ \dd + M^{-1} (\dd M) + M^{-1} A M \logd = \dd + B \logd \]
i.e.\ $ z \dd M = (M B - A M ) \dd z$.
In the Taylor expansion we must have
\[ j M^{j} = \sum_{k=0}^{j} \{ M^{k} B^{j-k} - A^{j-k} M^{k} \}
\,\,\,\,\,\,\,\,\,\,  \forall j\geq 0 .\]
So (\ref{one}) is fulfilled if
\begin{equation} \label{two}
  (j M^{j} + A^{0} M^{j} - M^{j} B^{0})-B^{j} = -A^{j} +
       \sum_{k=1}^{j-1} \{ M^{k} B^{j-k} - A^{j-k} M^{k} \} =:R^{j-1}
\end{equation}
for all $j \geq 1$.
Work by induction on $j \geq 1$.
For all $i,m=1, ..., l$, we need to satisfy, for $M^{0}=I $ ($A^{0}=B^{0}$),
 the equation on the block entries
\[ (j+A^{0}_{i,i})M_{i,m}^{j} - M_{i,m}^{j} A_{m,m}^{0} - B_{i,m}^{j} =
R^{j-1}_{i,m}. \]
By (\ref{three}), any eigenvalue $\lambda''$ of $(j+A_{i,i}^{0})$ satisfies
$[-\Real \lambda'']=-j+\psi^{i}$. So, $\lambda''$ is not an eigenvalue
of $A_{m,m}^{0}$ unless $\psi^{i}-j =\psi^{m}$.
Hence there is a solution $(M_{i,m}^{j},B_{i,m}^{j})$ with
$B_{i,m}^{j}=0$ if $\psi^{i} -j \neq \psi^{m}$.
So we find $M$ and $B$ as required.

\hfill

To see that $\sum_{0}^{\infty} M^{j} z^{j}$ is
absolutely convergent near $0$,
set  \[ c_{j}:= \| A^{j} \| + \| B^{j} \| \,\,\,\,\,\,\,\,\,\,
\forall j \geq 1 \,\,\,\,\,
\mbox{and} \,\,\,\,\,   c_{0}:=2 [ \, \| A^{0} \| \, ] +2 \in \Z \]
in operator norm. We can find $C > 1$ and $\varepsilon_{0} > 0$ such that
\[ c_{j} \varepsilon^{j} < C \,\,\,\,\,\,\,  \forall j\geq 0 \,\,\,\,\,\,\,
 \forall \,\,  0 \,\, \leq \varepsilon \leq \varepsilon_{0} \]
since $A$ and $B$ are absolutely convergent.
The equality (\ref{two}) implies
\[ (j-c_{0}) \| M^{j} \| \leq \sum_{k=0}^{j-1} \| M^{k} \| \,\, c_{j-k}. \]

Choose any $\delta \leq \varepsilon_{0}/2C$.
Then $(2C\delta)^{j-k} c_{j-k} \leq C$ for all
$j-k$ and so
\[ (j - c_{0})\| M^{j} \| (2C \delta)^{j} \leq C
\sum_{k=0}^{j-1} \| M^{k} \| (2C \delta) ^{k}. \]
Hence
\[ (j-c_{0}) \|M^{j} \| \delta^{j} \leq 2^{-j}
\sum_{0}^{j-1} \| M^{k} \| (2\delta)^{k}. \]

Let $D:= \sum_{k=0}^{c_{0}} \|M^{k} \| \delta^{k}$. Then we claim that
\[ \| M^{j} \| \delta^{j} \leq D \,\, 2^{c_{0}-j} \,\,\,\,\,\,\,\,\,\,
\forall j\geq c_{0} \]
which would finish (i).
The claim is clear for $j=c_{0}$ and for $j>c_{0}$ we use induction to find
\begin{eqnarray*}
\| M^{j} \| \delta^{j}  & \leq &
\frac{2^{-j}}{j-c_{0}} \left( \sum_{0}^{c_{0}} \|M^{k}\| \delta^{k} 2^{c_{0}}
        + \sum_{c_{0}+1}^{j-1} D 2^{c_{0}-k} 2^{k} \right) \\
  & \leq & 2^{c_{0}-j} D \frac{1+(j-1) -c_{0}}{j-c_{0}} .
\end{eqnarray*}

\hfill

{\bf (ii):}
By hypothesis,
\[ \nabla_{\theta'}=\dd +B(z) \frac{\dd z}{z}=\dd + u^{\Phi'}(-K'-\Phi')
  u^{-\Phi'} \frac{\dd u}{u}, \]
\[ \nabla_{\theta}=\dd +A(z) \frac{\dd z}{z}=\dd + z^{\Phi}(-K-\Phi)
  z^{-\Phi} \frac{\dd z}{z} \]
and $M=\sum_{0}^{\infty} M^{j}z^{j}:U \rightarrow \homo(\C^{r'},\C^{r})$
with $M \circ \nabla'_{\theta'}=\nabla_{\theta} \circ M$, i.e.\
(\ref{two}) above.

Let $M=\oplus_{i,m} M_{i,m}$ for $M_{i,m}:F'^{m} \rightarrow F^{i}$ where
$\C^{r'}=\oplus_{1}^{l'}F'^{m}$ and $\C^{r}=\oplus_{1}^{l} F^{i}$
according to $\Phi'$ and $\Phi$, respectively.
If $M$ is block-upper-triangular, i.e.\
$M_{i,m}=0$ if $\psi'^{m} > \psi^{i}$, then (ii) follows.

For $j=0$, (\ref{two}) gives $A^{0}M^{0}-M^{0}B^{0}=0$.
Clearly, $B^{0}$ and $A^{0}$ are block-diagonal, any eigenvalue
$\lambda'$ of $B_{m,m}^{0}$ satisfies $[-\Real \lambda']=\psi'^{m}$ and
any eigenvalue $\lambda$ of $A_{i,i}^{0}$ satisfies $[-\Real \lambda
 ] = \psi^{i}$. This implies
\begin{equation}
\label{four}
M_{i,m}^{0}=0 \mbox{\ \ if \ \ } \psi'^{m} \neq \psi^{i}.
\end{equation}
Since $A$ and $B$ are block-upper-triangular,
(\ref{two}) implies by induction
on $j \geq 0$ that $jM^{j} +A^{0}M^{j} -M^{j}B^{0}$ is
block-upper-triangular,
i.e. \[ (j+A_{i,i}^{0})M_{i,m}^{j} - M_{i,m}^{j}B^{0}_{m,m}=0 \]
if $\psi'^{m} > \psi^{i} > \psi^{i}-j$.
Hence, $M^{j}$ is block-upper-triangular.

\hfill

{\bf (iii):}
Assume $\tau$ is injective and $v' \in U \times \C^{r'}$
has $\phi'(v')=\psi'^{m}$. By (\ref{four}), there exists $i$ such that
$\psi^{i}=\psi'^{m}$ and $M^{0}_{i,m}$ has full rank.
Since $\phi$ and
$\phi'$ are invariant under parallel transport, $v'$ is, w.l.o.g.,
contained in the neighbourhood of $0$ where $M_{i,m}$ has full rank. Hence,
$M(v')$ has a component in $F^{i}$ and so, $\phi(M(v')) \leq \psi^{i}$.

Now assume $\tau$ is surjective and $v \in U \times \C^{r}$
with $\phi(v)=\psi^{i}$. Write $v=v^{1}+...+v^{i}$ according to the
decomposition of $\C^{r}$. By the dual of the above argument, we find
$v'=v'^{1}+...+v'^{i}$ with $v^{k}=M(v'^{k})$ and $\phi'(v'^{k}) =\psi^{k}$.
Hence, $\phi'(v')=\psi^{i}$.
\stopproof

\begin{defi}
Let $\theta$ be a normal trivialisation of $(E,\nabla)$. The integer
weights filtration of $(E,\nabla)$ is
\[ 0 \subset E^{1} \subset ... \subset E^{l} =E
\mbox{ \ \ where \ \ }  E^{m}:=\theta^{-1}(\oplus_{1}^{m} F^{k}). \]
\end{defi}

At first, the normal trivialisation $\theta$ and hence the filtration
of $E$ exists only over a small neighbourhood of the singularity. But since
each $E^{m}$ is invariant under $\nabla$, we can extend over all of $U$.
This filtration, together with $\Phi$, is equivalent to
$\phi \circ \theta:E \rightarrow \Z \cup \{ +\infty\}$.
It is independent of the choice of $\theta$ by Theorem \ref{nf}.

Let $\pi:\tilde{U}^{*} \rightarrow U^{*}$ be the universal covering
and write $\tz$ for the coordinate over $z$.
Let $\log \tz:\tilde{U}^{*} \rightarrow \C$ be a holomorphic function
such that $\log \tz \equiv \log z \mbox{ \ mod \ } (2 \pi i)$.
For $K \in \mbox{End}(\C^{r})$ let
$\tz^{K}:= \exp (K \log \tz)$; $\tz^{\Phi}=z^{\Phi}$.

\begin{lemm}
\showlabel{fs}
 For $\theta$ as in Definition \ref{nt},
   $ \nabla_{\theta} (z^{\Phi} \tz^{K} ) =0$
  on $\tilde{U}^{*}$,
i.e.\ $z^{\Phi} \tz^{K}$ is a fundamental system of flat sections.
 Hence, $\exp (2 \pi i K)$ is the monodromy around $0$.
\end{lemm}

\proof
$ (\dd + z^{\Phi} (-K -\Phi) z^{-\Phi} \frac{dz}{z} )(z^{\Phi} \tz^{K})
       =  z^{\Phi} ( \Phi + K - K - \Phi ) \tz^{K} \frac{dz}{z}. \,\, \Box $

\subsection{Correspondence between local logarithmic connections
and weighted flat bundles}

\begin{defi}[Simpson, Deligne]
A weighted flat bundle $(H,H^{m},\Phi)$ over $U^{*}$ consists
of a holomorphic, rank $r$ vector bundle $H \rightarrow U^{*}$ with a
holomorphic (i.e.\ flat and compatible) connection
\[ \nabla:\Omega^{0}(H) \rightarrow \Omega^{0}(H)
\otimes \Omega^{1}_{U^{*}},\]
 a filtration by proper subbundles
$ 0 \subset H^{1} \subset ... \subset H^{l}=H$,
invariant under $\nabla$,
and an $r \times r$ matrix with integer entries
$ \Phi= \mbox{diag }(\phi^{i})=\mbox{block-diag }(\psi^{m}I_{d^{m}}) $
where $\psi^{m} > \psi^{m+1}$ and $d^{m}=\rank (H^{m}/H^{m-1})$.
\end{defi}

The integer $\psi^{m}$ is called the {\em integer weight} of $H^{m}$.
The function
\function{\phi}{H}{\Z \cup \{ + \infty \}}{v}{\left\{
                    \begin{array}{lcl}
                      \psi^{m} & \mbox{if} & v \in H^{m} \setminus H^{m-1} \\
                      + \infty & \mbox{if} & v=0
                    \end{array}
                                            \right. }
is equivalent to the filtration and weights; $(H,\phi):=(H,H^{m},\Phi)$.
A {\em morphism} of weighted flat bundles
$\eta :(H',\phi') \rightarrow (H,\phi)$
is a map of flat bundles such that
\[ \phi(\eta(v')) \geq \phi'(v') \]
for all $v' \in H'$. Equivalently, $\eta(H'^{k}) \subseteq H^{m-1}$
if $\psi'^{k} > \psi^{m}$.
The morphism is called {\em injective} if it is as bundle map and satisfies
$ \phi(\eta(v'))=\phi'(v') $
for all $v' \in H'$.
It is called {\em surjective} if it is as bundle map
and if for all $v \in H$ there exists $v' \in \eta^{-1}(v) $
such that
$ \phi(v)=\phi'(v')$.
A sequence
\threehorss{(H',\phi')}{(H,\phi)}%
{(H'',\phi)}{\eta}{\xi}
is called {\em exact} if $\eta$ is injective, $\xi$ is surjective
and Im$\, \eta = \,\,$Ker$\, \xi$.
The {\em direct sum} $(H',\phi') \oplus (H'',\phi'')$ is defined by
$(H' \oplus H'',\phi)$ where $\phi(h' \oplus h''):= \min
(\phi'(h'),\phi''(h''))$.

\begin{defi}
Let $\calF^{-1}$ be the functor from the category of logarithmic connections
to the category of weighted flat bundles, given by restricting the
integer weights filtration and the morphisms to $U^{*}$.
\end{defi}

By Theorem \ref{nf},
$\calF ^{-1}$ sends injections, surjections and short exact
sequences to such.
To construct a functor $\calF$, inverse to $\calF^{-1}$,
consider a weighted flat bundle $(H,\phi)$.
Let $Y=(y_{1},...,y_{r})$ be a
fundamental system of multivalued flat sections,
 such that
$\langle y_{1},...,y_{d^{1}+...+d^{m}} \rangle = H^{m}$.

Let $\g$ be a loop in $U^{*}$ going once around $0$ anticlockwise
and write $\g ^{*}$ for the
induced action on $\tilu^{*}$; $\log (\tz \circ \gamma^{*})=
  (\log \tz)+2 \pi i$.
Then \[ Y \circ \g^{*} = Y G\]  for
 constant block-upper-triangular $G \in \Gl (r,\C)$.
Put  $K:= \normlog G$ which is also block-upper-triangular.
Since $\tz^{-K} \zmP$ is invertible over $\tilde{U}^{*}$,
$Y\tz^{-K} \zmP$ is a trivialisation of $H$ over $\tilde{U}^{*}$. It is
single valued, \cite[p 17]{ano}, since
\[ (Y \tz^{-K} \zmP)\circ \g^{*} = Y G G^{-1} \tz^{-K} \zmP.\]

\begin{defi}
\showlabel{f}
Let $\calF(H,\phi)$ be the
extension of $H$ over $U$,
whose stalk at $0$ is generated by the system of sections
\[ V(z)  := Y \tz^{-K} \zmP:U^{*} \rightarrow H \times \cdots \times H .\]
$\nabla$ becomes a singular  connection on the extension of $H$.
For a morphism
$ \eta: (H',\phi') \rightarrow (H,\phi),$
let $\calF(\eta): \calF(H',\phi') \rightarrow \calF (H,\phi)$
be the unique holomorphic extension of $\eta$.
\end{defi}

$\calF(H,0)$ is called the {\em canonical extension}, \cite{nas},
  of the flat bundle $H$.
One checks that
different choices of $Y$ and coordinate $z$ give
extensions which are isomorphic
via a map extending the identity of $H$.
Anosov \& Bolibruch construct extensions of $H$ by choosing $Y$ such that
$G$ is upper-triangular, but with $Y$
not requested to respect a fixed filtration.
One can (in addition)
choose $Y$ such that $G$ decomposes w.r.t.\ eigenvalues.
Assuming then that $G$ has only one eigenvalue, it is easy to see that
$\calF$ equals the extension-functor of Manin, \cite[p 94]{del},
and is a special case of Simpson's extension functor, \cite[p 738]{sim}.
They have the following lemma.

\begin{lemm}
$\calF (H,\phi)$
is a logarithmic connection and $\calF$ is inverse to $\calF^{-1}$ on
objects.
\end{lemm}

\proof
\hspace{0.5cm}  $ \nabla(V(z))  =$
\[ Y \dd (\tz^{-K} \zmP ) = Y \tz^{-K} (-K -\Phi) \zmP \logd
         = V(z) \zP (-K -\Phi) \zmP \logd. \]
So, $\nabla_{V} = \dd + \zP (-K -\Phi ) \zmP \frac{dz}{z}$
in the trivialisation given by the columns of $V$.
Combine this with Lemma \ref{fs}.
\stopproof

\begin{lemm}
\showlabel{gk}
Let $G,G' \in \Gl (r,\C)$ and put $K:=\normlog G$, $K':=\normlog G'$.
If  $C$ is an $r \times r'$-matrix
such that $ GC=CG'$
then $KC=CK'$ and hence
$ \tz^{K} C = C \tz^{K'}$.
If
$G G'=G' G$ then $KK'=K'K$.
$\Box$
\end{lemm}

\begin{lemm}[Simpson, Deligne]
\showlabel{he}
The holomorphic extension
$\calF(\eta)$ in Definition \ref{f} exists.
It commutes with the logarithmic connections.
If $\eta$ is injective (surjective) then so is $\calF(\eta) $.
$\calF$ sends short exact sequences to such and is inverse to $\calF^{-1}$
on morphisms.
\end{lemm}

\proof
Since $\eta$ maps flat sections to
flat sections we can find a unique constant
$r \times r'$-matrix $C$ such that
$ \eta \circ Y'=Y C$.
Because $\eta$ does not decrease weights,
 $\zP C z^{-\Phi'}$ is holomorphic over $U$.
Furthermore,
\[(\eta \circ Y') \circ \gamma^{*} = (YC) \circ \gamma^{*}
\mbox{\ \ so that \ \ }
 (\eta \circ Y')G'=YGC,\]
i.e.\ $CG'=GC$.
By Lemma \ref{gk} this implies
$ \tz^{K} C=C\tz^{K'}$.
We find
\[ \calF(\eta)(V')  =  (\eta \circ Y')\tz^{-K'} z^{-\Phi'}
                   =  V (\zP \tz^{K} C \tz^{-K'} z^{-\Phi'})
                  =  V(\zP C z^{-\Phi'}) \]
and $\calF(\eta)$ is holomorphic.
If $\eta$ is injective (surjective) then we can
choose $Y$ ($Y'$) such that $C$ is a permutation matrix
of full rank and $\zP C z^{-\Phi'} = C$.
The remainder of the statement follows from continuity.
\stopproof

\section{\sloppy
Global logarithmic connections and weighted flat bundles}

\subsection{Correspondence over $\Proj^{1}$}

\showlabel{no}
We extend the concepts to the Riemann sphere.
Let $X = \{ x_{1},...,x_{n} \} \subseteq  \Proj ^{1} $,
put $S:= \Proj^{1} - X$ and fix a base point $s \in S$.
For each $j=1,...,n$ choose a simply connected
neighbourhood $U_{j} \subseteq \Proj^{1}$
of $x_{j}$ containing $s$ but no other $x_{k}$'s and a coordinate $z_{j}$
centered at $x_{j}$.
($U_{j}$ is easier to handle than a small neighbourhood around $x_{j}$
and a path from $x_{j}$ to $s$.)
Let $\g_{j} \in \pi_{1}(U^{*}_{j},s)$ go once around $x_{j}$, anticlockwise,
and  $\tilu^{*}_{j}$ be the universal covering of $U^{*}_{j}$.

\hfill

A {\em logarithmic  connection} over $(\Proj^{1},X)$
 consists of a holomorphic bundle $E \rightarrow \Proj^{1}$
and a $\C$-linear map
$ \nabla:\Omega^{0}(E) \rightarrow \Omega^{0}(E) \otimes
\Omega^{1}( \log X)$
satisfying the Leibnitz rule; where $\Omega^{1}( \log X)=\Omega^{0}
(\calK_{\mathProj^{1}} \otimes [X])$.

A weighted flat bundle over $S$ is
a holomorphic flat bundle $H \rightarrow S$
together with filtrations by flat subbundles
$ 0 \subset H_{j}^{1} \subset ... \subset H_{j}^{l_{j}} = H|_{U^{*}_{j}}$
and integer weights
$ \Phi_{j}=\mbox{diag} (\phi^{i}_{j})=
\mbox{block-diag} (\psi_{j}^{m} I_{d_{j}^{m}})$
for each $j=1,...,n$.
Note that no compatibility is required for different $j$.
Write $\phi=(\phi_{1},...,\phi_{n})$ for the weight functions
$\phi_{j}:H|_{U_{j}^{*}} \rightarrow \Z \cup \{ +\infty \}$.

We have $\pi_{1}(S,s)= \langle \gamma_{1}, ..., \gamma_{n} \,\, | \,\,
 \gamma_{1} \cdot ... \cdot \gamma_{n} = 1 \rangle$
where $\gamma_{1} \cdot \gamma_{2}$ means travelling along
$\gamma_{1}$ first.
A weighted flat bundle corresponds to a conjugacy class of representations
$\chi:\pi_{1}(S,s) \rightarrow \Gl (r,\C)$ with,
for each $j=1,...,n$, a weighted filtration of
$\C^{r}$ invariant under $\chi(\gamma_{j})$.

\begin{coro}[\cite{sim}]
\showlabel{fc}
$\calF$ induces an
equivalence  between the category of logarithmic connections
over $\Proj^{1}$ and that of
weighted flat bundles over $S$.
It and its inverse preserve injections, surjections and short exact
sequences.
$\,\,\,\, \Box$
\end{coro}

\begin{defi}
\showlabel{fd}
\begin{description}
\item[(i)]
$ \deg (H,\phi) := \sum_{j} \{ \Tr \Phi_{j}
+ \Tr (\mbox{norm log} \,\,  \chi(\gamma_{j}))\} \in \Z $
\item[(ii)]
A system (in one of the considered categories) of rank $r$ and
degree $d$ is called stable if any
proper subsystem of
rank $r'$ and degree $d'$ satisfies $ d'/r' < d/r $.
For semi-stability allow '$ \leq $'.
The number $d/r$ is called the slope of the system.
\end{description}
\end{defi}

Observe that $\deg (H,\phi) = \sum \{ \Tr \Phi_{j}
+ \Real \Tr (\mbox{norm log} \,\,  \chi(\gamma_{j})) \}$
and it is an integer because
$\det \chi(\gamma_{1}) \cdot ... \cdot \det \chi(\gamma_{n})=1$.

Consider a weighted flat bundle $(H,\phi)$. Choose a basis
$Y_{s}$ of $H_{s}$ and denote its extensions by parallel transport
over $\tilde{U}_{j}^{*}$ by $Y(\tz_{j})$.
For each $j$, fix some $Z_{j} \in \Gl (r,\C)$ such that
$Y_{j}(\tz_{j}):=Y(\tz_{j})Z_{j}$ respects the filtration
of $H|_{U^{*}_{j}}$.
Let \[ Y_{j}(\tz_{j}) \circ \g_{j}^{*} = Y_{j}(\tz_{j})  G_{j}, \]
$K_{j} := \normlog G_{j}$.
Set $(E,\nabla)=\calF (H,\phi)$.
By the
Birkhoff-Grothendieck-Theorem (BGT), \cite{oss},
\cite{ano}, there is a system
\[ W:\Proj^{1} \rightarrow E \times ... \times E \]
of $r$ meromorphic sections such that $W|_{S}$ spans $H=E|_{S}$.
We have
\[ V_{j}(z_{j}):=Y_{j}(\tz_{j}) \tz_{j}^{-K_{j}} z_{j}^{-\Phi_{j}} =
W(z_{j}) Q_{j}(z_{j}),   \]
for some meromorphic $Q_{j}:U_{j} \rightarrow \Gl (r,\C)$ , holomorphic
on $U_{j}^{*}$.
Note that $W|_{U_{j}}$ spans $E|_{U_{j}}$ if and
only if $Q_{j}$ is holomorphic at
$x_{j}$.

\begin{prop}[{\cite[p 32]{lev}}, {\cite[p 754]{sim}}]
\showlabel{de}
$\calF$ preserves degrees and so \\ (semi-) stability.
\end{prop}

\proof
Let $\omega$ be the connection matrix of
$\nabla|_{H}$ w.r.t.\ the trivialisation
$W|_{H}$.
$\Tr \omega$ is a single valued holomorphic one-form
on $S$. By the Residue-Theorem,
$0=\sum_{1}^{n} \mbox{Res} _{x_{j}} (\Tr \omega ).$
Since
$ \nabla(W|_{U_{j}^{*}})= Y_{j}(\tz_{j}) \dd (\tz_{j}^{-K_{j}}
z_{j}^{-\Phi_{j}} Q_{j}(z_{j})^{-1}),$
\begin{eqnarray*}
 \omega|_{U_{j}^{*}}& =& Q_{j} (z_{j}) z_{j}^{\Phi_{j}} \tz_{j}^{K_{j}}
\dd ( \tz_{j}^{-K_{j}} z_{j}^{-\Phi_{j}} Q_{j}(z_{j})^{-1}) \\
                & = & -Q_{j} [ z_{j}^{\Phi_{j}} K_{j} z_{j}^{-\Phi_{j}}
               + \Phi_{j} ] Q_{j}^{-1} \frac{ \dd z_{j}}{z_{j}}
                         + Q_{j} \dd (Q_{j}^{-1}).
\end{eqnarray*}
Let $k_{j}$ be the order of vanishing of $\det Q_{j}$ at $x_{j}$.
Then
\[ \Tr \omega|_{U_{j}^{*}}
 = - (\Tr K_{j} + \Tr \Phi_{j} + k_{j} ) \frac{\dd z_{j}}{z_{j}} + \alpha \]
for a holomorphic 1-form $\alpha$. Hence
$ -\sum_{j} k_{j} = \sum_{j} \Tr K_{j}  + \Tr \Phi _{j} .$
The right-hand-side is the degree of the weighted
flat bundle, while the left-hand-side is the sum of the orders
of vanishing of $\det W$,
the degree of $E$.
\stopproof

This result holds, in fact, over any Riemann surface.
Note, that it implies Lemma 5.2.2 in \cite{ano}.
Also,
if $(E,\nabla) \rightarrow \Proj^{1}$ is logarithmic at $X$
and has residues $\rho_{j}:E_{x_{j}} \rightarrow E_{x_{j}}$
then $-\sum_{1}^{n} \Tr \rho_{j}$ equals the degree of $E$.
This is because
in a normal trivialisation, using the notation of Definition \ref{nt},
$\rho_{j}=(z_{j}^{\Phi_{j}}(-K_{j}-\Phi_{j}) z_{j}^{-\Phi_{j}})(0)$.
Hence, $-\sum \Tr \rho_{j}=\sum (\Tr K_{j}+\Tr \Phi_{j})=\deg \calF^{-1}
 (E,\nabla )$.

\subsection{The splitting type of $E \rightarrow \Proj^{1}$}

By the BGT, any holomorphic bundle
$E \rightarrow \Proj^{1}$
has the form
$ E \cong \calO (c_{1}) \oplus ... \oplus \calO (c_{r}) $
for unique integers $c_{1} \geq ... \geq c_{r}$.
We call $C:= \diag (c_{i})$ the {\em splitting type} of $E$ or $(E,\nabla)$.
If $C=c_{1} I_{r}$ we say $E$ has {\em constant}
splitting type.
Recall that $\HH^{0}(\calO(c))$ is zero if $c <0$ and equal to $\C$ if $c=0$.
As $\calF$ preserves subsystems and degrees, this implies
Theorem 5.2.2 in \cite{ano}.

\begin{theo}
\showlabel{ne}
If $(E,\nabla)$ is semi-stable, $n \geq 2$ and $C=\mbox{diag}\,(c_{i})$
the splitting type, then
$ (0\leq ) \,\, c_{i} - c_{i+1} \leq n-2 $ for all $i=1,...,r-1.$
\end{theo}

\proof
Fix a splitting
$ E = \calO(c_{1}) \oplus ... \oplus \calO(c_{r}) . $
Suppose there exists an $i \in \{1,...,r-1 \}$
such that $c_{i} - c_{i+1} > n-2$.
Let $z$ be a coordinate centred at $s \in S$.
For all $k \in \{1,...,i \}$ we can find
sections $v_{k}$ of $\calO (c_{k})$  vanishing to order $c_{k}$ at $s$,
i.e.\ $v_{k} z^{-c_{k}}$ spans $\calO (c_{k})$ near $s$.
For each $m \in \{ i+1,...,r\}$ we consider the natural projection
\[ \pi_{m} : E \otimes \calK_{\mathProj^{1}} \otimes [X] \rightarrow
\calO (c_{m}) \otimes \calK_{\mathProj^{1}} \otimes [X] \]
and obtain  sections
$\pi_{m} \circ \nabla (v_{k}) : \Proj^{1} \rightarrow \calO (c_{m})
\otimes \calK_{\mathProj^{1}} \otimes [X].$
Near $s$ we have
\begin{eqnarray*}
 \pi_{m} \circ \nabla (v_{k}) & =
         & \pi_{m} \circ \nabla (v_{k} z^{-c_{k}} z^{c_{k}} ) \\
                              & =
 & \pi_{m} ( \nabla ( v_{k} z^{-c_{k}}) z^{c_{k}} + (v_{k} z^{-c_{k}})
                                    \dd (z^{c_{k}} )) \\
                & = & \pi_{m} ( \nabla (v_{k} z^{-c_{k}} ) z^{c_{k}} )
\end{eqnarray*}
since $k \neq m$.
Either $\pi_{m} \circ \nabla (v_{k})$ is
identically zero or of order at least $c_{k}$.
The latter is equivalent to
\[ c_{k} \leq \deg (\calO (c_{m}) \otimes
\calK_{\mathProj^{1}} \otimes [X]) = c_{m} +n - 2 . \]
But from $k < i<m$ we see that
$ c_{k} \geq c_{i} > c_{i+1} +n -2 \geq c_{m} +n -2 .$
Therefore, $\pi_{m} \circ \nabla (v_{k})$ is
identically zero for all $k<i<m$. So, $\nabla$
preserves $\calO (c_{1}) \oplus ... \oplus \calO(c_{i})$.
Semistability implies then that $c_{1} =c_{2}=...=c_{r}$.
\stopproof

This theorem easily  extends to
logarithmic connections with parabolic structure at the
singularities, i.e.\ to filtered regular $D_{S}$-modules, c.f. \cite{sim}.
Note, any logarithmic connection with irreducible monodromy is semi-stable,
even stable.

\begin{coro}[Bolibruch, {\cite[p 84]{ano}}]
\showlabel{bo}
If the monodromy of $(E,\nabla)$ is irreducible ($n \geq 2$) then
$ \sum_{i=1}^{r} c_{1} - c_{i} \leq (n-2)r(r-1)/2$. $\Box$
\end{coro}

\begin{lemm}
\showlabel{ct}
If $\calF (H,H_{j}^{m},\Phi_{j})$ has splitting type
$C$ and $\Phi'_{j}=\Phi_{j}+\lambda_{j}I_{r}$ then
$\calF(H,H_{j}^{m},\Phi'_{j}) $
has splitting type $C+(\sum_{1}^{n} \lambda_{j})I_{r}$.
$\Box$
\end{lemm}

\begin{lemm}
\showlabel{cc}
Let
$ \threehorbb{(H',\phi')}{(H,\phi)}{(H'',\phi'')}$
be a short exact sequence
and assume that two of them have the same slope.
Then all three have the same slope and
$\calF (H,\phi)$ has constant splitting type
$C=c I_{r}$ if and only if $\calF (H',\phi')$ and $\calF (H'',\phi'')$
have constant splitting types
$C'=cI_{r'}$ and $C''=cI_{r''}$, respectively. $\Box$
\end{lemm}

The proofs of these two lemmas are straightforward.
\section{The Riemann-Hilbert-Problem}

\subsection{Commutative  and semi-simple representations}

\begin{lemm}[\cite{bol}, {\cite[p 76]{ano}}]
\showlabel{co}
If $\chi:\pi_{1}(S) \rightarrow \Gl\, (r,\C)$
factors through $\HH_{1}(S)$  then it is the monodromy
of a Fuchsian system.
\end{lemm}

\proof
Since the $G_{j}=\chi(\gamma_{j})$ commute, each $G_{j}$ preserves
$\ker (G_{k}-\mu I)^{t}$ of each $G_{k}$ and for each $t$.
Assume then that each $G_{j}$ has
only one eigenvalue $\rho_{j}$ and is upper-triangular.
Let
$\mu_{j}:=\mbox{norm log}\,\, \rho_{j}$
be the only eigenvalue of $K_{j}:=\mbox{norm log}\,\, G_{j}$
and $\xi:=\sum_{1}^{n} \mu_{j} \in \Z$.
By Lemma \ref{gk} the $K_{j}$'s commute and
$\exp(\xi \cdot I_{r} -\sum K_{j})=G_{1} \cdot ... \cdot G_{n} =I_{r}$.
Since $\xi \cdot
I_{r}-\sum K_{j}$ has only the eigenvalue $0$, it is the normalised
logarithm of $I_{r}$, i.e. $0$.
A short calculation then shows that
\[ \nabla:= \dd + \left[ \frac{\xi}{z-x_{1}} -
\sum_{j=1}^{n} \frac{K_{j}}{z-x_{j}} \right] \dd z \]
is smooth at infinity.
Over each $\tilde{U}_{k}^{*}$ we set
$ Y:= (z-x_{1})^{-\xi} \prod_{1}^{n}
        \widetilde{(z-x_{j})}^{K_{j}}$ and find
$\nabla (Y) = 0$ and $Y \circ \g_{j} = Y G_{j}$.
\stopproof

\begin{prop}[{\cite[p 80]{ano}}]
\showlabel{bq}
If $(E,\nabla)=\calF(H,\phi)$ has splitting type $C$
and $k \in \{1,...,n \}$ is fixed
then there exists a permutation $P$ and meromorphic
$W:\Proj^{1} \rightarrow E \times ... \times E$ such that
\begin{description}
\item[(i)] $W$ is holomorphic except at $x_{k}$
   and spans $E$ away from $x_{k}$,
\item[(ii)] $Wz_{k}^{-C}$ spans $E$ near $x_{k}$ and
\item[(iii)]
 $ Q_{k} (z_{k}) = \hat{Q}_{k}(z_{k}) z_{k}^{-C} P =
\hat{Q}_{k}(z_{k}) P z_{k}^{-P^{-1}CP} $
for some invertible $\hat{Q}_{k}$.
\end{description}
\end{prop}

We give a proof different from that in \cite{ano}.

\proof
Let $W'=(w'_{1},...,w'_{r})$ where $w_{i}':\Proj^{1} \rightarrow
\calO (c_{i})$ vanishes to order $c_{i}$ at $s$ and the $\calO(c_{i})$'s
decompose $E$.
Define $Q_{k}'$  by $V_{k}=W' Q'_{k}$ near $x_{k}$, $\det Q_{k} \neq 0$.

Claim: For each permutation $P$ such that all
bottom-right minors of $Q_{k}'(0) P^{-1}$ are non-singular,
there exists a $W$ as in the proposition.
(The existence of such $P$ follows by induction from the description of
the determinant of a matrix in terms of co-rank one minors.)

We may assume that $z_{k}(s)=\infty$.
Suppose there exists
$ b=((b_{i,j})):\Proj^{1} \rightarrow \Gl (r,\C)$
such that
\begin{equation}
\label{i}
 b_{i,j} = \left\{ \begin{array}{cc}
           \sum_{0}^{c_{i}-c_{j}-1} b_{i,j}^{p} z_{k}^{p}  & i<j \\
            1                                           & i=j \\
            0                                         & i>j,
                    \end{array}
          \right.
\end{equation}
\begin{equation}
\label{ii}
 z_{k}^{c_{i} -c_{m}} \,\, | \,\, (b \, Q'_{k} P^{-1})_{i,m} \,\, \ \ \
\forall \,\,\, i<m.
\end{equation}
Then $W=W' b z_{k}^{C}$
spans $E|_{\mathProj^{1}-\{ x_{k} \}}$, $Wz_{k}^{-C}$
spans $E$ near $x_{k}$ and
$W' Q_{k}'=V_{k}=W Q_{k}$ implies
\[ Q_{k}=z_{k}^{-C} b Q'_{k} = z_{k}^{-C} (b Q'_{k} P^{-1})P
=\hat{Q}_{k} z_{k}^{-C} P \]
for some invertible $\hat{Q}_{k}$ and we would be done.
To find $b$ as in (\ref{i}) satisfying (\ref{ii})
we need to solve a system of linear equations.
With $Q'_{k}P^{-1}=((q_{j,m}))$, condition (\ref{ii}) is equivalent to
\[ z_{k}^{c_{i}-c_{m}} \ | \
\sum_{j=i+1}^{r} b_{i,j} q_{j,m} + q_{i,m} \]
for all $1 \leq i < m \leq r$.
Writing
$q_{j,m} = \sum_{0}^{\infty} q_{j,m}^{p} z_{k}^{p}$,
(\ref{ii}) becomes
\[ \sum_{j=i+1}^{r}
\sum _{t=0}^{c_{i}-c_{j}-1} b_{i,j}^{t} q_{j,m}^{p-t} + q_{i,m}^{p} =0 \]
for all $1\leq i < m \leq r$ and $0 \leq p <c_{i}-c_{m}$.
We define
\[\alpha(t):=\min \{ j \in \{i+1,...,r \}
\ | \ t \leq c_{i} - c_{j} -1 \}. \]
Then (\ref{ii}) is equal to
\[ \sum_{t=0}^{p} \sum_{j=\alpha(t)}^{r} b_{i,j}^{t} q_{j,m}^{p-t}
+ q_{i,m}^{p}=0 \]
for all $1 \leq i < m \leq r$ and
$0 \leq p < c_{i}-c_{m}$. Note that $\alpha(p) \leq m $.
To find the $b_{i,j}^{t}$'s we argue  one row of $b$ at a time, i.e.\
fix $i$.
Assume that $b_{i,j}^{t}$ is known by induction for all $t<p$.
Then we have to fullfil
\[ \sum_{j=\alpha(p)}^{r} b_{i,j}^{p} q_{j,m}^{0} = \mbox{known term} \]
for $m=\alpha(p),...,r$.
This system has a solution $(b_{i,\alpha(p)}^{p},...,b_{i,r}^{p})$
since the matrix $((q_{j,m}^{0}))_{\alpha(p) \leq j,m \leq r}$
is a right-bottom minor of
$Q'_{k}(0)P^{-1}$.
\stopproof

\begin{coro}[Plemelj]
\showlabel{ss}
If there exists $k \in \{ 1,..., n\}$ such that
$\chi (\g_{k})$ is semi-simple (i.e.\ diagonalizable)
then $\chi$ is the monodromy of a Fuchsian system.
\end{coro}

\proof
If $\chi (\g_{k})$ is semi-simple we can split
\[H|_{U_{k}^{*}} = H_{k,1} \oplus ... \oplus H_{k,r} \]
into flat line bundles.
Let $(E,\nabla)$ be the canonical extension of $H$
(i.e.\ $\Phi_{j}=0$ for all $j$)
and choose $Z_{k}$
(see subsection \ref{no}) such that the i-th section in $Y_{k}(\tz_{k})=
Y(\tz_{k}) Z_{k}$ spans $H_{k,i}$.
Choose $W$ as in  Proposition \ref{bq}.
So $Q_{k}=\hat{Q}_{k}z_{k}^{-C} P$
for invertible $\hat{Q}_{k}$ and permutation $P$, where
$C$ is the splitting type of $E$.

Let $P'$ be the permutation with $P'_{i,j}=0$ if
$i + j \neq n$ and $P'_{i,n-i}=1$ for $i=1,...,n$.
Consider the filtration of $H|_{U_{k}^{*}}$, induced by the
sections in $Y_{k}(\tz_{k})P^{-1} (P')^{-1}$. As $\chi(\gamma_{k})$ is
semi-simple it will respect this filtration.
Let $\Phi'_{k}:=-(P')^{-1} C P'$, which is diagonal
with non-increasing entries,
 and $\Phi'_{j} =0$ for $j \neq k$.
Put $(E',\nabla'):=\calF (H,\Phi'_{j})$.
Then $W$  trivializes $E'|_{\mathProj^{1}-x_{k}}$.
Furthermore, $Q_{k}'$ is invertible and hence,
$W$ spans $E'$ globally, since
\begin{eqnarray*}
W Q'_{k} & = & Y_{k}(\tz_{k}) P^{-1} (P')^{-1} \tz_{k}^{-K_{k}'}
               z_{k}^{-\Phi_{k}'} = Y_{k}(\tz_{k}) \tz_{k}^{-K_{k}} P^{-1}
               (P')^{-1} z_{k}^{-\Phi'_{k}} \\
         & = & W Q_{k} P^{-1} z_{k}^{C} (P')^{-1} = W \hat{Q}_{k}
               (P')^{-1}. \,\, \Box
\end{eqnarray*}

\subsection{The rank two case}

For a representation $\chi:\pi_{1}(S,s) \rightarrow \Gl \, (2,\C)$ with
canonical extension $(E^{0},\nabla^{0})=\calF (H,0)$
of splitting type $C^{0}=\mbox{ diag }(c_{1}^{0},c_{2}^{0})$, Bolibruch calls
$c_{1}^{0}-c_{2}^{0}$ the {\em weight} of the canonical extension,
\cite[p 102]{ano}.

\begin{theo}
\showlabel{tw}
\begin{description}
\item[(i) (Dekkers)]
Any rank two representation $\chi$ is the holo\-nomy
of a Fuchsian system.
\item[(ii) (Bolibruch, {\cite[p 137]{ano}})]
$ c_{1}^{0}-c_{2}^{0}= \min_{\phi} \sum_{1}^{n}(\phi_{j}^{1} -
    \phi_{j}^{2})$
where $\phi$ runs over all integer weight functions on $H$ such that
$\calF (H,\phi)$ is Fuchsian.
\end{description}
\end{theo}

\proof
Assume $(E,\nabla)=\calF (H,\phi)$ is Fuchsian.
The identity of $H$ extends to a meromorphic map
$\calO(c_{1}^{0}) \oplus \calO(c_{2}^{0})=E^{0} \rightarrow E$.
By definition of $\calF$ (or the proof of Lemma \ref{he}), the non-zero
map $\calO(c_{1}^{0}) \rightarrow E$ is of order greater or equal
$\phi_{j}^{2}$ at $x_{j}$. Hence,
$c_{1}^{0} \leq - \sum_{1}^{n} \phi_{j}^{2}$.
By Proposition \ref{de},
$c_{1}^{0}+c_{2}^{0}+\sum (\phi_{j}^{1} + \phi_{j}^{2})=0$ and hence
$c_{1}^{0}-c_{2}^{0} \leq \sum_{1}^{n} (\phi_{j}^{1} - \phi_{j}^{2})$.
This proves (ii) in one direction.

If there exists $k \in \{1,...,n\}$ such that $\chi(\gamma_{k})$
has two eigenvalues then we are done by the proof of Corollary \ref{ss}.
Otherwise, $(E^{0},\nabla^{0})=\calF (H,0)$ is semi-stable
and we can argue much as in Theorem \ref{he}.
If $v_{i}:\Proj^{1} \rightarrow \calO(c_{i}^{0})$ vanishes
to order $c_{i}^{0}$
at  $s$ and $\pi_{2}:E^{0} \rightarrow \calO (c_{2}^{0})$,
then \[ \pi_{2} \circ \nabla^{0}(v_{1}):\Proj^{1} \rightarrow
 \calO (c_{2}^{0} ) \otimes \calK_{\mathProj^{1}} \otimes [X] \]
either vanishes identically, in which case $c_{1}^{0}=c_{2}^{0}$ and
we are done,
or it is of order at least $c_{1}^{0}$
at $s$. If it also vanishes at each $x_{j}$ then
$n+c_{1}^{0} \leq c_{2}^{0}+n-2$.
So, assume $(\pi_{2} \circ \nabla^{0} (v_{1})(x_{k}) \neq 0$.
If $W=(v_{1},v_{2})$, the formula for
the connection-matrix $\omega|_{U^{*}_{k}}$ in the
proof of Proposition \ref{de} implies that $(Q_{k})_{1,2}(x_{k}) \neq 0$
($\det Q_{k} (x_{k}) \neq 0$ by choice of $W$).

Now apply the Claim at the beginning of the proof of Proposition
\ref{bq} with $W'=(v_{1},v_{2})$ and $P$ the non-trivial rank two
permutation. Then $\calF (H,0,...,0,\Phi_{k},0,...,0)$ will be Fuchsian
for $\Phi_{k}=-P^{-1} C^{0} P$.
We have completed the proof of (ii) and also proved (i).
\stopproof

\subsection{The semi-stable case and implications}

\begin{theo}
\showlabel{ta}
Let $H \rightarrow S$ be a flat bundle.
If we can find filtrations $H_{j}^{m}$ and integer weights $\Phi_{j}$
such that
\begin{description}
\item[(a)] $(H,H_{j}^{m},\Phi_{j})$ is semi-stable and
\item[(b)] there exists $k \in \{1,...,n\}$ with
        $\mbox{\rm rank} (H_{k}^{i+1} / H_{k}^{i})=1$
       and \[ \phi_{k}^{i} - \phi_{k}^{i+1} \geq (r-1)(n-2)
             \,\,\,  \forall \,\,\, i=1,...,r-1 \]
\end{description}
then we can find $\Phi_{k}'$ such that
$ (E',\nabla'):=\calF (H,H_{j}^{m},\Phi_{1},...,\Phi_{k}',...,\Phi_{n}) $
is Fuchsian.
\end{theo}

\proof
Let $C$ be the splitting type of $(E,\nabla) = \calF (H,H_{j}^{m},\Phi_{j})$.
Theorem \ref{ne} implies
\[ c_{1} - c_{r} \leq (n-2)(r-1) \leq \phi_{k}^{i} - \phi_{k}^{i+1} \]
for all $i=1,...,r-1$.
Fix $W$ as in  Proposition \ref{bq} and let, in that notation,
$\Phi_{k}':= \Phi_{k}-P^{-1}CP$, which will have non-increasing entries.
$W$ spans $E'$ off $x_{k}$ and
\[WQ_{k}z_{k}^{\Phi_{k}} \tz^{K_{k}}_{k} = Y_{k}(\tz_{k}) =
             W Q'_{k} z_{k}^{\Phi_{k}'} \tz_{k}^{K_{k}} \]
implies
\[ Q_{k}'=Q_{k} z_{k}^{P^{-1}CP} =
        \hat{Q}_{k} z_{k}^{-C} P z_{k}^{P^{-1}CP} = \hat{Q}_{k}P. \]
So, $Q_{k}'$ is invertible and
$W$  a global trivialisation of $E'$.
\stopproof

\begin{prop}
\showlabel{ca}
Suppose there exists $k \in \{1,...,n \}$ and
$h \in H_{s}$  such that $h$ is an eigenvector of $\chi(\g_{k})$ but
a cyclic vector of the $\pi_{1}(S,s)$-module $H_{s}$ (i.e.\
$ \langle \,\, (\mbox{Im} \,\, \chi) (h) \,\, \rangle = H_{s}$).
Let $N_{1},...,N_{n}$ be any integers.

Then we can find
filtrations $H_{j}^{m}$ and weights $\Phi_{j}$ such that
$\calF (H,H_{j}^{m},\Phi_{j})$ is Fuchsian.
Moreover, we can arrange that
\begin{description}
\item[(1)] $\phi_{j}^{i} \geq N_{j}$ for all $j \neq k$, $i=1,...,r$ and
\item[(2)] $\phi_{k}^{1}=\phi(h)\geq N_{k}$.
\end{description}
Hence, there are infinitely many Fuchsian systems with monodromy $\chi$.
\end{prop}

\proof
For each $j\neq k$ choose filtrations $H_{j}^{m}$ and weights $\Phi_{j}$
such that (1) is satisfied.
Also choose a filtration $H_{k}^{m}$ and weights $\Phi_{k}$
such that $H_{k}^{1}=\langle h \rangle$,
hypothesis (b) of Theorem \ref{ta}
is satisfied and $\phi_{k}^{1} \geq N_{k}+(r-1)(n-2)$.
These conditions remain satisfied if we increase $\phi_{k}^{1}$ or decrease
$\phi_{k}^{r}$. Doing so we can assume that
$ \deg (H,H_{j}^{m},\Phi_{j}) =0 $
and since no proper flat subbundle of $H$ contains $h$
we can also assume that
$(H,H_{j}^{m},\Phi_{j})$ is semi-stable.

Apply Theorem \ref{ta} to find
$ \sum_{1}^{r} c_{i}=0 $ and $ c_{i}-c_{i+1} \leq (n-2), $
implying $|c_{i}| \leq (n-2)(r-1)$ for $i=1,...,r$.
Hence, the first entry of $\Phi'_{k}=\Phi_{k}-P^{-1}CP$
is greater or equal to
$N_{k}+(r-1)(n-2)-(r-1)(n-2)=N_{k}$ and we are done.
\stopproof

\begin{coro}[Bolibruch, {\cite[p 84]{bol}}, {\cite[p 83]{ano}}; Kostov,
\cite{kos}]
\showlabel{tb}
Any \\ irreducible flat $H \rightarrow S$
is the restriction of a Fuchsian system.
$\Box$
\end{coro}

\begin{coro}[{\cite[p 114]{ano}}]
Any
$\chi: \pi_{1}(S,s) \rightarrow \Gl(r,\C)$ is the subrepresentation
of the monodromy of some Fuchsian system of double the rank.
\end{coro}

\proof
By Lemma \ref{co} we may assume that $n \geq 3$ and $r \geq 2$.
Let $G_{j}:= \chi (\g_{j})$ for all $j=1,...,n$.
By Corollary \ref{ss} we can assume that in canonical basis
$e_{i}=(0,...,0,1,0,...,0)^{t}$ we have the equality of vector spaces
$ \langle e_{r},G_{1}e_{r} \rangle = \langle e_{r-1},e_{r} \rangle$.
We define $G'_{j}$ as follows
\[ G_{1}':= \twomatrix{G_{1}}{M_{1}}{0}{I} \ \ \ \
        M_{1}:= \left( \begin{array}{cccccc}
                        1 & & &  & & \\
                         & \cdot & 0 & & & \\
                        & 0 & \cdot & & & \\
                         & & & 1 & & \\
                          & & & & 0 & 0 \\
                         & & & & 1 & 0
                          \end{array}
                \right) \]
\[ G_{2}':= \twomatrix{G_{2}}{0}{0}{M_{2}} \ \ \ \
         M_{2}:= \left( \begin{array}{ccccc}
                       1 & 1 & & & \\
                       & \cdot &  \cdot & 0 & \\
                         & & \cdot & \cdot  & \\
                         & 0 & & \cdot & 1 \\
                       & & & & 1
                       \end{array}
                  \right)
          \]
\[ G_{3}':=
         \twomatrix{G_{3}}{-G_{2}^{-1} G_{1}^{-1} M_{1}}{0}{M_{2}^{-1}}
\mbox{\  \ and \ \ }
 G'_{j}:= \twomatrix{G_{j}}{0}{0}{I} \]
for all $j \geq 4$.
One checks that
$G_{1}' \cdot ... \cdot G_{n}' = I_{2r}$ and so defines
a representation $\chi':\pi_{1}(S,s) \rightarrow%
\Gl (2 r,\C)$.
Furthermore, $e_{2r}$ is an eigenvector of $G_{1}'=\chi'(\gamma_{1})$
and
\[ \langle e_{2r},G_{2}' e_{2r},...,(G_{2}')^{r-1} e_{2r} \rangle =
 \langle e_{r+1},...,e_{2r} \rangle ,\]
\[ G_{1}' \langle e_{r+1},...,e_{2r} \rangle
\oplus \langle e_{r+1},...,e_{2r} \rangle
 \supset \langle e_{1},...,e_{r-2},e_{r},...,e_{2r} \rangle \]
and
$ \langle e_{r},G_{1}' e_{r} \rangle = \langle e_{r-1},e_{r} \rangle .$
Apply Proposition \ref{ca} with $h:=e_{2r}$.
\stopproof

\subsection{Reducible representations}

Let $H_{j}:=H|_{U_{j}^{*}}$.
Part (i) of the following is due to Bolibruch,
\cite[Cor. 5.4.1]{ano}, \cite[Thm. 3.8]{bol},
 while (ii) is new and will
be used to give an alternative proof of Bolibruch's answer
to the RHP
in rank three, Theorem \ref{dr} here.

\begin{theo}
\showlabel{re}
Let \threehorbb{H'}{H}{H''} be a short exact sequence of flat bundles
(without weights)
and assume there exist filtrations and weights such that
 $(E'',\nabla'')=\calF (H'',(H'')_{j}^{m},\Phi''_{j})$
is Fuchsian.
Suppose at least one of the following
conditions holds for some $k \in \{1,...,n\}$.
\begin{description}
\item[(i)]
$\threehorbb{H'_{k}}{H_{k}}{H''_{k}}$ splits and
there exist filtrations and weights such that
$(E',\nabla')=\calF(H',(H')_{j}^{m},\Phi'_{j})$ is Fuchsian.
\item[(ii)] There exist splittings
   $ H'_{k}=H^{(3)} \oplus H^{(0)}$ and
          $H_{k}=H^{(3)} \oplus H^{(4)}$
where $\langle h \rangle =H^{(0)}=H'_{k} \cap H^{(4)}$
and $h$ is a cyclic vector of the $\pi_{1}(S,s)$ module $H'_{s}$.
\end{description}
Then $H$ is the restriction of a Fuchsian system.
\end{theo}

\proof
{\bf (i):}
By Lemma \ref{ct} we can assume
that for all $j \neq k$ the smallest diagonal entry in $\Phi'_{j}$
is greater then the largest one in $\Phi''_{j}$.
Using \[ \threehorbb{H'_{j}}{H_{j}}{H''_{j}} \]
for such $j$ we can therefore induce filtrations and weights
on $H_{j}$ to make this local sequence
a short exact one of weighted flat bundle.

Let $\alpha:H''_{k} \rightarrow H_{k}$ be a splitting right
inverse of $H_{k} \rightarrow H''_{k}$
and put the obvious weighted filtration on $\Image \alpha$.
Then use
$H_{k}=H'_{k} \oplus \Image \alpha$ to give $H_{k}$ the direct
sum weighted filtration.
$(H'_{k},\phi'_{k}) \rightarrow (H_{k},\phi_{k})$
becomes an injection. Since
$(H_{k},\phi_{k}) \rightarrow (H''_{k},\phi''_{k})$ is the composition
$ \threehorbs{H_{k}}{\Image \alpha}{H''_{k}}{\alpha^{-1}} $,
it is a surjection.
Apply Lemma \ref{cc} to finish this case.

\hfill

{\bf (ii):}
For $j=1,...,n$ let $N_{j}$ be the greatest diagonal entry in $\Phi''_{j}$.
Then construct $(H',(H')_{j}^{m},\Phi'_{j})$ as in Proposition \ref{ca}
so that $\calF(H',(H')_{j}^{m},\Phi'_{j})$ is Fuchsian.
If $j\neq k$ we induce weights on $H_{j}$ as in (i).

Use the exact sequence of flat bundles
$ \threehorbb{H^{(0)}}{H^{(4)}}{H''_{k}}$
to induce weights $\phi^{(4)}$ on $H^{(4)}$.
Induce weights on $H_{k}$ using
$ H_{k}=H^{(3)} \oplus H^{(4)}$.
Then, $(H_{k},\phi_{k}) \rightarrow (H''_{k},\phi''_{k})$
is given by the composition
\[ \threehorbb{(H_{k},\phi_{k})}{(H^{(4)},\phi^{(4)})}%
{(H''_{k},\phi''_{k})} \]
of two surjections and hence is one itself.
For $h'=h^{(3)} + h^{(0)} \in H'_{k}$ we have
\[ \phi'_{k}(h')=\min (\phi'_{k}(h^{(3)}),\phi'_{k}(h^{(0)})) =
  \phi(h') \]
since $H^{(0)}$ is the highest weight subspace in the
filtration of $H'_{k}$.
Hence $(H'_{k},\phi'_{k}) \rightarrow (H_{k},\phi_{k})$
is an injection and we finish with
Lemma \ref{cc}.
\stopproof

\begin{prop}[{\cite[p 83]{bol}}, {\cite[p 100]{ano}}]
\showlabel{bt}
If $H$ is reducible, the holonomy $\chi(\gamma_{j})$
has only one Jordan-block,
for each $j=1,...,n$, and $\calF(H,\phi)=(E,\nabla)$
has constant splitting type, then $\Phi_{j}=\phi_{j}^{1} \cdot I_{r}$
for all $j$.
\end{prop}

We give a more conceptual proof.

\proof
If $G_{j}$ has only one Jordan block,
there exists a canonical
full flag of subsystems of $H_{j}$.
If $H'$ is a proper subsystem of $H$
then it must contain the subbundles in the local filtrations
of rank equal to rank $H'$.
If some $\Phi_{k}$
has non equal diagonal entries
then
\[ \mbox{ slope}\,\, (H',\phi')>\, \mbox{slope}\,\, (H,\phi).\]
But $(E,\nabla)$, and hence $(H,\phi)$,
is semi-stable -- a contradiction.
\stopproof

Let $H$ be as in the previous proposition,
$\rho_{j}$ the only eigenvalue
of $\chi(\gamma_{j})$ and $\mu_{j}=\,\, \mbox{norm log}\,\, \rho_{j}$.
If $(E,\nabla)=\calF(H,\phi)$ is Fuchsian then
$\deg (H,\phi)=0$.
Since $r \, | \, \sum \Tr \Phi_{j}$ we find
$r \, | \, \sum \Tr (\mbox{norm log } \chi (\gamma_{j}))$, i.e.\
$r \, | \, r \sum \mu_{j}$.
So, $\sum \mu_{j}$ must be an integer.
Bolibruch uses this to give an example of a representation
with $r= 4$ and $n=3$ which can not be the monodromy of any Fuchsian system
with three singularities, \cite[p 91]{bol}, \cite[p 105]{ano}.

\subsection{Parabolic representations and the rank three case}

Let $\para(r,\C)$ be the group of invertible upper-triangular
 $r \times r$-matrices.

\begin{theo}
\showlabel{p}
Let $\chi:\pi_{1}(S,s) \rightarrow \para(r,\C)$ be a representation with
\[ G_{j}=\chi(\gamma_{j})=\uppermatrix{\rho^{1}_{j}}{\rho_{j}^{r}} \]
for $j=1,...,n$.
Let
 $ \mu_{j}^{i}:=\mbox{\rm norm log} (\rho_{j}^{i}) $
 and
 $\Lambda^{i}:=-\sum_{j=1}^{n} \Real \mu_{j}^{i} \in \Z_{\leq 0}.$
Assume we can find $((\phi_{j}^{i}))$
such that
\begin{description}
\item[(a)] $\phi_{j}^{i} \geq \phi_{j}^{k}$ if ($i \leq k$ and
$\rho_{j}^{i}=\rho_{j}^{k}$) and
\item[(b)] $\Lambda^{i}=\sum_{j=1}^{n} \phi_{j}^{i}$ for all $i=1,...,r$.
\end{description}
Then $\chi$ is the monodromy of a Fuchsian system $(E,\nabla)$
and the
integer weights of $(E,\nabla)$ equal $((\phi^{i}_{j}))$ as sets.
\end{theo}

\proof
The flat bundle $H$, associated to $\chi$,
has a global natural filtration.
We work by induction on the rank $r$
and extend the claim of the theorem by the fact that
the integer weights function $\phi_{j}$ on $H_{j}$,
which we construct, is the direct sum of its
restrictions to the generalised eigenspaces
of $G_{j}$, acting on $H_{j}$.

For $r=1$ we let
$\phi_{j}:H_{j}-\{ 0\} \rightarrow \Z \,\,\,$
have single value $\phi_{j}^{1}$.
This implies  $\deg (H,H_{j}^{m},\Phi_{j})=\Lambda^{1}-\Lambda^{1}=0$
and hence $\calF (H,H_{j}^{m},\Phi_{j})$ is
Fuchsian.

For $r \geq 2$
write $H' \subseteq H$, $H'_{j} \subseteq H_{j}$ for
the rank $(r-1)$ subbundles.
We are given, by induction,
\[ \phi_{j}':H_{j}' \rightarrow %
\{ \phi_{j}^{1},...,\phi_{j}^{r-1} \} \cup \{ +\infty \} \]
with the above described property.
Consider one $j \in\{1,...,n\}$ at a time.
Let $A \subseteq H_{j}$ ($A' \subseteq H_{j}'$) be the generalised
eigenspace of $G_{j}$ ($G'_{j}$)
of the eigenvalue $\rho_{j}^{r}$.
Also let $B' \subseteq H_{j}'$ be the direct
sum of the other generalised eigenspaces of $G'_{j}$.
Then, the extended induction hypothesis implies
$ \phi'_{j}(h_{A'}+h_{B'})=
\min (\phi_{j}' (h_{A'}),\phi'_{j})(h_{B'})).$
Put
\function{\phi_{A}}{A}{\Z \cup \{ + \infty\}}{h}{\left\{ \begin{array}{ccc}
                            \phi'_{j}(h) & \mbox{if} & h \in A' \\
                            \phi_{j}^{r} & \mbox{if} & h \in A-A'.
                                                      \end{array}
                                                 \right.}
We give  $H_{j}=A \oplus B'$ the direct sum of the weighted filtrations.
By construction, $(H'_{j},\phi_{j}') \rightarrow (H_{j},\phi_{j})$
is an injection.
We give $H_{j}/H'_{j}$ the integer weight $\phi_{j}^{r}$.
Then $h=h_{A}+h_{B'} \in H_{j}$ maps to zero under
$\alpha:H_{j} \rightarrow H_{j}/H'_{j}$
unless $h_{A} \in A-A'$ in which case
\[ \phi(h)=\min (\phi_{A}(h_{A}),\phi'_{j} (h_{B'})) \leq
\phi_{A}(h_{A}) = \phi_{j}^{r}. \]
$\alpha$ is surjective since
$\phi(h_{A})=\phi_{j}^{r}$ for any $h_{A}\in A-A'$.
We have constructed a short exact sequence
\[ \threehorbb{(H',\phi')}{(H,\phi)}{(H/H',\phi^{H'/H})} \]
where $H'$ is the restriction of a Fuchsian system
by induction and $H/H'$ is so similar to the rank one case.
Apply Lemma \ref{cc} to finish.
\stopproof

The following result is due to Bolibruch when $r=3$, \cite[p 133]{ano}.
He has a counter example for
$r=7$, $n=4$,
\cite[p 106]{ano}.

\begin{coro}
\showlabel{th}
For $\chi: \pi_{1}(S,s) \rightarrow \para(r,\C)$
 and $r\in \{1,2,3,4\}$ there exists a Fuchsian system with monodromy $\chi$.
\end{coro}

\proof
We want to find
$((\phi_{j}^{i}))$
satisfying (b) of Theorem \ref{p} and \\
{\bf (a)':} $\phi_{j}^{i} =\phi_{j}^{k}$ if $\rho_{j}^{i}=\rho_{j}^{k}$.

Claim : If there is $m \in \{1,...,n\}$ and $k \in \{1,...,r\}$ such that
$\rho_{m}^{k} \neq \rho_{m}^{t}$ for all $t \neq k$ then the problem to find
$((\phi_{j}^{i}))$, satisfying (b) and (a)', reduces to rank $(r-1)$.

To see this just find $((\phi_{j}^{i}))_{i\neq k}$,
choose $(\phi_{j}^{k})_{j \neq m}$ if they are not
fixed by (a)' already and calculate $\phi_{m}^{k}$ using (b).

The corollary is trivial for $r=1$. For $r=2,3$,
we can either use the claim or
have $\Lambda^{1}=...=\Lambda^{r}$ and finish easily.
For $r=4$, if we can not use the claim there are two cases.
Either,
for each $j=1,...,n$,
($\rho_{j}^{1}=\rho_{j}^{3}$ and $\rho_{j}^{2}=\rho_{j}^{4}$)
or ($\rho_{j}^{1}=\rho_{j}^{4}$
and $\rho_{j}^{2}=\rho_{j}^{3}$).
Hence, $\Lambda^{1}+\Lambda^{2} -\Lambda^{3}=\Lambda^{4}$.
Solve the system consisting of (a)' and (b)
for $((\phi_{j}^{i}))_{j=1,..,n;i=1,2,3}$
and get a solution for the  rank four problem by setting
$\phi_{j}^{4}:= \phi_{j}^{1}+\phi_{j}^{2}-\phi_{j}^{3}$ for $j=1,...,n$.

Or, $\rho_{m}^{1}=\rho_{m}^{2}\neq \rho_{m}^{3}=\rho_{m}^{4}$ for some
$m \in \{1,...,n\}$. Solve two rank two problems, i.e.\
find $((\phi_{j}^{i}))_{i=1,2}$ satisfying (a)' and (b)
and find $((\phi_{j}^{i}))_{i=3,4}$ satisfying (a)' and (b).
Increasing all entries of $((\phi_{j}^{i}))_{i=1,2;j\neq m}$ by
a sufficiently large integer $N$ and decreasing $\phi_{m}^{1}$ and
$\phi_{m}^{2}$ by $(n-1)N$, we can satisfy (a) and (b) for the rank four
problem.
\stopproof

\begin{theo}[{\cite[p 90]{bol}}, {\cite[p 133]{ano}}]
\showlabel{dr}
A rank three representation $\chi:\pi_{1}(S,s) \rightarrow \Gl \, (3,\C)$
is the monodromy of a Fuchsian system
if and only if one or more of the following holds.
\begin{description}
\item[(a)] $\chi$ is irreducible.
\item[(b)] Some $\chi(\gamma_{k})=G_{k}$ has more than one Jordan block.
\item[(c)] The canonical extension of $H=H(\chi)$
 has constant splitting type.
\end{description}
\end{theo}
The part of the proof which is left we  do differently from Bolibruch.

\proof
By Corollary \ref{tb}, Corollary \ref{ss}, Lemma \ref{ct} and
Proposition \ref{bt} we are left to
show that if $\chi(\gamma_{k})$ has
two Jordan blocks for some $k$ then $\chi$ is the monodromy
of a logarithmic connection on $\Proj^{1} \times \C^{3}$.

Let $h_{1} \in \C^{3}$ ($h_{2} \in \C^{3}$) be the
eigenvector corresponding to the
rank one (rank two) Jordan block of $\chi(\gamma_{k})$.
Consider the $\pi_{1}(S,s)$-submodules
\[ F_{1}:=\langle (\Image \chi)(h_{1}) \rangle \,\,\,\,\,
 \mbox{and} \,\,\,\,\,
   F_{2}:= \langle (\Image \chi )(h_{2}) \rangle \]
of $\C^{3}$.
If $\rank F_{1}=1$ we are in case (i) of Theorem \ref{re}  because of
the positive solvability of the RHP in rank two.
If $\rank F_{1}=3$ then we are in the case
of Proposition \ref{ca}.
So assume $\rank F_{1}=2$ and hence $h_{2} \in F_{2} \subseteq F_{1}$
and $\rank F_{2} <3$.
If $\rank F_{2}=1$ we are in the case of Corollay
\ref{th} and if $\rank F_{2}=2$
we are in case (ii) of Theorem \ref{re}
with $H'=F_{1}$, $H^{(0)}=\langle h_{2} \rangle$ and
$H^{(3)}=\langle h_{1} \rangle$.
\stopproof


\end{document}